\documentclass[reqno,11pt]{amsart}
 
\usepackage{color}

\newtheorem{theorem}{Theorem}[section]
\newtheorem{lemma}[theorem]{Lemma}

\theoremstyle{definition}
\newtheorem{assumption}[theorem]{Assumption}

\theoremstyle{remark}
\newtheorem{remark}[theorem]{Remark}

\makeatletter
\def\dashint{\operatorname%
{\,\,\text{\bf--}\kern-.98em\DOTSI\intop\ilimits@\!\!}}
\makeatother

\newcommand\bR{\mathbb{R}}

\newcommand\cB{\mathcal{B}}
\newcommand\cF{\mathcal{F}}

\newcommand\cM{\mathcal{M}}
\newcommand\cN{\mathcal{N}}
\newcommand\cP{\mathcal{P}}

\newcommand{\osc}{{\rm osc}\,}
\newcommand{\loc}{{\rm loc}\,}

 \newcommand{\mysection}[1]{\section{#1}
 \setcounter{equation}{0}}
 
\newcommand{\nlimsup}{\operatornamewithlimits{\overline{lim}}}

\begin{document}

\title%
{On strong solutions of It\^o's equations
with a$\,\in W^{1}_{d}$ and
   b$\,\in L_{d}$}
\author{N.V. Krylov}

\email{nkrylov@umn.edu}
\address{127 Vincent Hall, University of Minnesota, Minneapolis, MN, 55455}
 
 

\begin{abstract} 
We consider It\^o uniformly nondegenerate equations
with time independent coefficients, the diffusion
coefficient in $W^{1}_{d,\loc}$, and the drift in
$L_{d}$. We prove the unique strong solvability
for any starting point and prove that as a function
of the starting point the solutions are H\"older
continuous with any exponent $<1$. We also prove that
if we are given a sequence of coefficients 
converging in an appropriate sense to the original ones,
then the solutions of approximating equations
converge to the solution of the original one.
\end{abstract}

\maketitle

\mysection{Introduction}
                                                  \label{section 3.11.1}
Let $\bR^{d}$ be a $d-$dimensional Euclidean space of points
$x=(x^{1},...,x^{d})$ with $d\geq3$. Let $(\Omega,\cF,P)$ be a 
complete probability space,
let $\{\cF_{t}\}$ be an increasing filtration of 
$\sigma$-fields $\cF_{t}\subset \cF$, that are complete.
Let  
$w_{t}$ be a $d_{1}$-dimensional Wiener process relative to
$\{\cF_{t}\}$, where $d_{1}\geq d$.

Assume that on $\bR^{d}$ we are given $\bR^{d}$-valued Borel
functions $b,\sigma^{k}=(\sigma^{ik})$, $k=1,...,d_{1}$. We are going
to fix $x_{0}\in\bR^{d}$ and  investigate
the equation
\begin{equation}
                                         \label{6.15.2}
 x_{t}=x_{0}+\int_{0}^{t}\sigma^{k}(x_{s})\,dw^{k}_{s}+
\int_{0}^{t}b(x_{s})\,ds,
\end{equation}
where and everywhere below the summation 
over repeated indices is understood.

We are interested in the so-called strong solutions, that is
solutions such that, for each $t\geq0$, $x_{t}$ is $\cF^{w}_{t}$-measurable,
where $\cF^{w}_{t}$ is the completion of $\sigma(w_{s}:s\leq t)$. We present
sufficient conditions for the equation to have a strong solution
and also for the solution
to be unique (strong uniqueness).
A very reach literature
on the weak uniqueness problem for \eqref{6.15.2}     is beyond
  the scope of this article.

After the classical work by It\^o showing that there exists
a unique strong solution of \eqref{6.15.2} if $\sigma^{k}$ and $b$
are Lipschitz continuous (may also depend on time and $\omega$),
many efforts were applied to relax these Lipschitz conditions.
In   case $d=d_{1}=1$ T. Yamada and S. Watanabe \cite{YW_71} relaxed
the Lipschitz condition on $\sigma$ to the H\"older $(1/2)$-condition
(and even slightly weaker condition) and kept $b$ Lipschitz
(slightly less restrictive). Much attention was payed to equations
with continuous coefficients
satisfying the so-called monotonicity conditions
(see, for instance, \cite{Kr_84} and the references therein).

T. Yamada and S. Watanabe \cite{YW_71} also put forward
a very strong theorem, basically, saying that strong uniqueness
implies the existence of strong
solutions. Unlike the present paper, 
the majority of papers on the subject after that time
are using their theorem.
 S. Nakao (\cite{Na_72}) proved
the strong solvability in time homogeneous case
 if $d=d_{1}=1$ and $\sigma$ is bounded away from zero and infinity
 and is
locally of bounded variation. He also assumed that $b$ is bounded,
but from his arguments it is clear that the summability of $|b|$
suffices. In this respect his result shows that
our results are also true if $d=1$. However, the 
general case that $d=2$ is quite open.

A. Veretennikov  seems to be the first who in \cite{Ve_80} 
not only proved the existence of strong solutions
in the time inhomogeneous multidimensional case  when $b$ is bounded,
but also considered the case of $\sigma^{k}$ in Sobolev class,
namely,   $\sigma^{k}_{x}\in L_{2d,\loc}$. He used A. Zvonkin's method (see
\cite{Zv_74}) of transforming the equation in such a way that the
drift term disappears.
X. Zhang in \cite{Zh_11} considered time inhomogeneous equations
under some conditions
which for the time homogeneous case (our case) become 
$ \sigma^{k}_{x}, b\in L_{p}$ with $p>d$. For more detailed
information on the time inhomogeneous case we refer the reader
to \cite{Zh_11},  \cite{Zh_20}, and the references therein.  

Even the case when the $\sigma^{k}$'s are constant and the process is nondegenerate
attracted
very much attention especially in the time inhomogeneous
setting. We discuss some of the results in the particular case of
$ b$ independent of $t$.
M. R\"ockner and the author in \cite{KR_05} proved,
among other things, the existence of strong solutions
when $b\in L_{p}$ with $p>d$.
If $b$ is bounded  A. Shaposhnikov (\cite{Sh_14}, \cite{Sh_17})
proved the so called path-by-path uniqueness, which, basically,
means that for almost any trajectory $w_{t}$ there is only one
solution (adapted or not). This result was already announced
by A. Davie before with a very entangled proof which left many doubtful.

In a fundamental work by L.~Beck, F.~Flandoli, M.~Gubinelli, and M.~Maurelli
(\cite{BFGM_19}) the authors investigate such equations from
 the points of view of It\^o stochastic equations, stochastic transport 
equations, and stochastic continuity equations. Their article
contains an enormous amount of information and a vast references list.
We compare only those of their results which
have counterparts in the present article.
In what concerns our situation they require ($\sigma^{k}$ constant
and the process is nondegenerate) $b\in L_{p,\loc}$ with $p>d$
or $p=d$ but $\|b\|_{L_{p}}$ to be sufficiently small
and they prove strong solvability and strong uniqueness
(actually, path-by-path-uniqueness which is stronger) only for almost
all starting points $x$. We assume that $\sigma^{k}_{x},b\in L_{d}$
and for uniformly nondegenerate   and bounded $\sigma^{k}$
prove that, for any $x$, equation \eqref{6.15.2}
has a unique strong solution.

Our approach is absolutely different from
all articles mentioned above and all articles
which one can find in their references.
 We do not use
Yamada-Watanabe theorem or transformations of the noise.
Instead, our method is based on
an analytic criterion for the existence
of strong solutions which first appeared in \cite{VK_76}.

Simple examples of equations for which we prove the existence
of  unique strong solutions are
$$
dx_{t}=\big(2+I_{x_{t}\ne0}\zeta(x_{t})\sin(\ln|\ln |x_{t}|)\big) \,dw_{t},\quad
dx_{t}=dw_{t}+\zeta(x_{t})( |x_{t}|\ln |x_{t}|)^{-1}  \,dt,
$$
where $\zeta$ is any smooth function vanishing for $|x|>1/2$
satisfying $|\zeta|\leq 1$
and $l$ is any vector in $\bR^{d}$. Observe that in the first equation
the diffusion coefficient is discontinuous at the origin.

We conclude the introduction by some notation.
We set $u_{x}=Du$ to be the gradient of $u$, $u_{xx}$ to be
the matrix of its second-order derivatives,
$$
D_{x^{i}}=D_{i}u=u_{x^{i}}=\frac{\partial}{\partial x^{i}}u,\quad
u_{x^{i}\eta^{j}}=D_{x^{i}\eta^{j}}u=
D_{x^{i}}D_{\eta^{j}}u,
$$
$$
 \partial_{t}u=\frac{\partial}{\partial t}u\quad
u_{(\xi)}=\xi^{i}u_{x^{i}}.
$$
If $\sigma(x)=(\sigma^{i}(x))$ is vector-valued
(column-vector), by $\sigma_{x}$ we mean the matrix
whose $ij$th element is $ \sigma_{x^{j}}^{i}$. If $c$ is a matrix
(in particular, vector), we set $|c|^{2}=\text{tr}cc^{*}$
($=\text{tr}c\bar c^{*}$ if $c$ is complex-valued).

For $p\in[1,\infty)$ by $L_{p}$ we mean the space 
of Borel (perhaps complex- vector- or matrix-valued)
 functions on $\bR^{d}$ with finite norm given by
$$
\|f\|_{L_{p}}^{p}=\int_{\bR^{d}}|f(x)|^{p}\,dx.
$$
By $W^{2}_{p}$ we mean the space 
of Borel functions $u$ on $\bR^{d}$ whose Sobolev derivatives 
$u_{x}$ and $u_{xx}$ exist and $u,u_{x},u_{xx}\in L_{p}$.
The  norm in $W^{2}_{p}$ is given by
$$
\|u\|_{W^{2}_{p}}=\|u_{xx}\|_{L_{p}}+\|u \|_{L_{p}}.
$$ 
Similarly  $W^{1}_{p}$ is defined. As usual, we write $f\in L_{p,\loc}$
if $f\zeta\in L_{p}$ for any $\zeta\in C^{\infty}_{0}$ ($=C^{\infty}_{0}
(\bR^{d})$). Similarly $W^{\cdot}_{\cdot,\loc}$ are defined.

If a Borel $\Gamma\subset \bR^{d}$, by $|\Gamma|$ we mean its Lebesgue
measure.  Finally,
$$
B_{R}(x)=\{y\in\bR^{d}:|x-y|<R\},\quad B_{R} =B_{R}(0).
$$
 
\mysection{Main results}

Set $a^{ij}=\sigma^{ik}\sigma^{jk}$, $a=(a^{ij})$. Fix numbers $\delta\in(0,1)$
and $\|b\|,\| \sigma^{k}_{x}\|\in(0,\infty)$.

\begin{assumption}
                                                  \label{assumption 3.1.1}

We have
\begin{equation}     
                            \label{7.18.1}
\delta^{-1}|\lambda|^{2}\geq
a^{ij}(x) \lambda^{i}\lambda^{j}\geq\delta|\lambda|^{2}
\end{equation}
for all $\lambda,x\in\bR^{d}$.
Also   
 $$
\|b\|_{ L_{d} }\leq\|b\| .
$$
\end{assumption}

\begin{assumption}
                                                  \label{assumption 3.1.2} 
For any $k$ we have $\sigma^{k}\in W^{1}_{d,\loc}$ and
$$
\|\sigma^{k}_{x}\|_{L_{d}}\leq \| \sigma^{k}_{x}\|.
$$
 
\end{assumption}

\begin{theorem}
                                                     \label{theorem 6.29.1}
Under the above assumptions, for any $x_{0}\in \bR^{d}$, equation \eqref{6.15.2}
has a strong solution $x_{t}  $. If $y_{t}$ is also a solution of
\eqref{6.15.2}, then with probability one $x_{t}=y_{t}$ for all $t$.
\end{theorem}

\begin{theorem}
                                                     \label{theorem 6.29.3}
Under the above assumptions
suppose that we are also given sequences $\sigma^{k}(n)$,
$b(n)$, $n=1,2...$, $k=1,...,d_{1}$,
 of functions having the same meaning as $\sigma^{k}$,
$b$ and satisfying Assumptions \ref{assumption 3.1.1}
and \ref{assumption 3.1.2} with the same $\delta$, $\|b\|$ and $\| \sigma^{k}_{x}\|$.
Assume that $b(n)\to b$ and $\sigma^{k}_{x}(n)\to\sigma^{k}_{x}$
in $L_{d}$ as $n\to\infty$
and we are given  a sequence $x(n)\to x_{0}$. Finally, let $\sigma^{k} (n)\to\sigma^{k} $ (a.e.)
as $n\to\infty$. Then for any $m,T\in(0,\infty)$  
$$
\lim_{n\to\infty}E\sup_{t\leq T}| x_{t}(n,x(n)) - x_{t}|^{m}=0,
$$
where $x_{t}(n,x(n))$ are the solutions of \eqref{6.15.2}
in which $x_{0}$, $\sigma^{k}$, and $b$ are replaced by $x(n)$, $\sigma^{k}(n)$, and $b(n)$,
respectively.

\end{theorem}

\begin{theorem}
                                                     \label{theorem 6.29.2}
Under the above assumptions, 
there is a function $x_{t}(x)=x_{t}(\omega,x)$ which for  $x=x_{0}$
is a solution of \eqref{6.15.2} and for each $\alpha<1$  and $\omega$
is $\alpha$-H\"older continuous with respect to $x$ and $(\alpha/2)$-H\"older 
continuous with respect to $t$ on each set $[0,T]\times \bar B_{R}$,
$T,R\in(0,\infty)$.
 
\end{theorem}

\begin{remark}
                                                 \label{remark 7.11.1}
The main emphasis of the article is to treat the case that
$b\in L_{d}$. It is known (see, for instance, \cite{BF_61}
\cite{Kr_19_1}) that, even if $d_{1}=d$ and $(\sigma^{k})$
is a unit matrix, there are cases when $b\in L_{d-\varepsilon}$
for any $\varepsilon\in(0,1)$, but not for $\varepsilon=0$,
and there are no solutions of \eqref{6.15.2}.

 However, our results are new also if $b$
is bounded or $b\equiv0$. In that case the arguments  
are not so technically involved
 and   allow any $d\geq 1$ rather than $d\geq 3$.
In Remark \ref{remark 7.10.1} we show an example with $b\equiv0$
and $\sigma^{k}\in L_{d-\varepsilon}$
for any $\varepsilon\in(0,1)$, but not for $\varepsilon=0$,
when there are no strong solutions. In this regard Assumption
\ref{assumption 3.1.2} seems to be optimal.
\end{remark}
 
The rest of the article is organized as follows.
As we mentioned above our main tool is an 
analytic criterion for the existence of strong solutions.
To derive it we develop  necessary facts from the theory
of semigroups generated by elliptic operators in Section \ref{section 3.11.2}.
Then in Section \ref{section 6.18.1} we relate the semigroup
from Section \ref{section 3.11.2} to the semigroup of the corresponding
Markov diffusion process. In Section \ref{section 7.3.1}
we derive our analytic criterion. Section \ref{section 7.3.2}
is devoted to some estimates of the series involved in the criterion
when $\sigma^{k}$ and $b$  are smooth. In Sections \ref{section 7.3.4},
\ref{section 7.3.5},  and \ref{section 7.3.9} we prove Theorems \ref{theorem 6.29.1},
\ref{theorem 6.29.3}, and \ref{theorem 6.29.2}, respectively.

\mysection{An analytic semigroup}
                                                  \label{section 3.11.2}
In this section Assumption \ref{assumption 3.1.1}
is supposed to be satisfied but Assumption
\ref{assumption 3.1.2} is replaced with a weaker Assumption
\ref{assumption 2.20.1} which comes after some discussion.

 Introduce the uniformly  elliptic operators
  $$
Lu( x)=(1/2) a^{ij}( x) u_{x^{i}x^{j}} ( x)+
b^{i}( x) u_{x^{i}}( x),
$$
$$
L_{0}u( x)=(1/2) a^{ij}( x) u_{x^{i}x^{j}} ( x)
$$
acting on functions given on $\bR^{d}$. 

Denote
$$
\osc (a,B_{\rho}(x))=
 |B_{\rho}|^{-2}
 \int_{y,z\in B_{\rho}(x)}|a( y)-a( z)|\,dydz,
$$
$$
a^{\# }_{r}=\sup_{ x \in\bR^{d }}\sup_{\rho< r}
\osc (a,B_{\rho}(x)) .
$$
 
Here is a consequence of Theorem 1 of \cite{DK_11}. We are dealing with 
complex-valued functions and
denote $\bR^{d+1}=\{(x^{0},x^{1},...,x^{d}):x^{k}\in\bR\}$.
\begin{lemma}
                          \label{lemma 2.19.2}
For any $p\in(1,\infty)$ and $\varepsilon\in(0,1]$  there exists
$\theta_{0}=\theta_{0}(d,\delta,\varepsilon,p)$ such that, if there is $r_{0}>0$
for which
$a^{\#}_{r_{0}}\leq\theta_{0}$, then there exist
  $\lambda_{0}\geq 1, N_{0}$,
depending only on $d,\delta,\varepsilon,p,r_{0}$, 
  such that, for any $u\in W^{2}_{p} $
and $\lambda\geq \lambda_{0}$,
$$
\sum_{r,s=0}^{d}\|D_{rs}u\|_{L_{p}(\bR^{d+1}) }  +\lambda\| u\|_{L_{p}(\bR^{d+1})  }
$$
\begin{equation}
                                 \label{6.4.1}
\leq N_{0} \|L_{0}u\pm \varepsilon iD^{2}_{0}u-(1\pm\varepsilon i)\lambda
u\|_{L_{p}(\bR^{d+1})  }.
\end{equation}
 \end{lemma}

Proof. As is easy to see, Theorem 1 of \cite{DK_11} is applicable
to the operator $Mu:=(1\pm\varepsilon i)^{-1}\big(L_{0}u\pm\varepsilon
iD^{2}_{0}u\big)$ and it yields an estimate for
$$
u\in  W^{1,2}_{p}((-\infty,0)\times\bR^{d+1})=\{u\in L_{p}
((-\infty,0)\times\bR^{d+1}):\partial_{t}u,
$$
$$
u_{x},u_{xx}
\in L_{p}
((-\infty,0)\times\bR^{d+1})\}
$$
 similar to \eqref{6.4.1}  where $\bR^{d+1}$
is replaced with $(-\infty,0)\times\bR^{d+1}$ and 
$L_{0}u\pm\varepsilon iD^{2}_{0}u-(1\pm\varepsilon i)\lambda u$ is replaced with
$Mu-\partial u/
\partial t-\lambda u$. By substituting in this estimate $u(x)e^{t}$,
we get \eqref{6.4.1} and the lemma is proved.

\begin{remark}
                                                \label{remark 6.3.1}
Without introducing the new coordinate, Theorem 1 of \cite{DK_11}
implies that in Lemma \ref{lemma 2.19.2} one can replace
\eqref{6.4.1} with
\begin{equation}
                                 \label{6.4.2}
\sum_{r,s=1}^{d}\|D_{rs}u\|_{L_{p} }+ \lambda\| u\|_{L_{p} }
\leq N_{0} \|L_{0}u - \lambda
u\|_{L_{p} }
\end{equation}
valid for any $u\in W^{2}_{p} $ and
$\lambda\geq \lambda_{0}(d,\delta, p,r_{0})$ with
$N_{0}=N_{0} (d,\delta, p,r_{0})$ as long as 
$a^{\#}_{r_{0}}\leq\theta_{0}(d,\delta, p)$.
Since $\|L_{0}u - \lambda
u\|_{L_{p} }\leq \|L_{0}u - \mu
u\|_{L_{p} }+|\lambda-\mu|\,\|u\|_{L_{p} }$,
estimate \eqref{6.4.2} easily implies that in the same situation
\begin{equation}
                                 \label{6.4.3}
\sum_{r,s=1}^{d}\|D_{rs}u\|_{L_{p} } +|\mu|\,\| u\|_{L_{p} }
\leq 2N_{0} \|L_{0}u - \mu
u\|_{L_{p} }
\end{equation}
as long as $|\mu|\geq \lambda_{0}$ and $|\Im \mu| \leq 2\varepsilon_{0}\Re \mu$,
where $\varepsilon_{0}=\varepsilon_{0}(d,\delta, p,r_{0})>0$.
\end{remark}

\begin{lemma}
                          \label{lemma 6.3.1}
For any $p\in(1,\infty)$   there exists $\theta_{0}=\theta_{0}(d,\delta,p)$
such that, if there is $r_{0}>0$ for which $a^{\#}_{r_{0}}\leq\theta_{0}$,
then there exist
  $\lambda_{0}\geq 1, N_{0}$,
depending only on $d,\delta,p,r_{0}$, 
  such that, for any $u\in W^{2}_{p} $
and complex $\lambda$ such that $\Re\lambda\geq \lambda_{0}$,
\begin{equation}
                                 \label{6.4.5}
\sum_{r,s=1}^{d}\|D_{rs}u\|_{L_{p} } +|\lambda|\| u\|_{L_{p} )}
\leq N_{0} \|L_{0}u -\lambda u\|_{L_{p} }.
\end{equation}
 
 \end{lemma}

Proof. We use an idea from \cite{Ag_62}.
Take a nonnegative $\zeta\in C^{\infty}_{0}(\bR)$ such that $\zeta^{p}$
has unit integral,
$u\in W^{2}_{p} $,
and plug into \eqref{6.4.1} the function $u(x)e^{i\mu x_{0}}\zeta(x_{0})$
and $\varepsilon=\varepsilon_{0}$.
Then we get for $\lambda\geq\lambda_{0}$ and $\mu\in\bR$ that
$$
\sum_{r,s=1}^{d}\|D_{rs}u\|_{L_{p} } + (\lambda+\mu^{2}) \|
u\|_{L_{p} }-N(1+|\mu|)\|
u\|_{L_{p} }
$$
\begin{equation}
                                 \label{6.4.4}
\leq N  \|L_{0}u -\big[(1\pm\varepsilon_{0} i)\lambda \pm\varepsilon_{0}
i\mu^{2}\big]u\|_{L_{p}} +N(1+|\mu|)\|
u\|_{L_{p}}.
\end{equation}
Now   take $\hat\lambda$ such that $\Re 
 \hat\lambda\geq  \lambda_{0}$. If $|\Im \hat\lambda|
 \leq 2\varepsilon_{0}\Re
\hat\lambda$, we have \eqref{6.4.5} with $\hat\lambda$
in place of $\lambda$ thanks to \eqref{6.4.3}.

If $\Im \hat\lambda
 \geq 2\varepsilon_{0}\Re
\hat\lambda$ set $\lambda=\Re \hat\lambda$, $\varepsilon_{0}\mu^{2}=
\Im \hat\lambda-\varepsilon_{0}\lambda$. Then
$$
|\hat\lambda|^{2}\leq ((2\varepsilon_{0})^{-2}+1)(\Im \hat\lambda)^{2},\quad
 \mu^{2}\leq \varepsilon_{0}^{-1}\Im \hat\lambda\leq
\varepsilon_{0}^{-1}| \hat\lambda|,\quad
\lambda+\mu^{2}=\varepsilon_{0}^{-1}\Im \hat\lambda
$$
and \eqref{6.4.4} with upper signs yields
$$
\sum_{r,s=1}^{d}\|D_{rs}u\|_{L_{p} } + |\hat\lambda|\, \|
u\|_{L_{p}} 
\leq N  \|L_{0}u -\hat\lambda u\|_{L_{p}} +N(1+|\hat\lambda|^{1/2})\|
u\|_{L_{p}}.
$$
By increasing $\lambda_{0}$ we absorb the last term on the right
into the left-hand side for $\Re 
 \hat\lambda\geq  \lambda_{0}$ and we come to
\eqref{6.4.5} with $\hat\lambda$ in place of $\lambda$
if $\Im \hat \lambda\geq0$. The case of $\Im \hat \lambda\leq0$
is treated by using
 \eqref{6.4.4} with lower signs.
The lemma is proved.

The argument in the second part of Remark \ref{remark 6.3.1}
also allows us to deduce from Lemma \ref{lemma 6.3.1}
the following.

\begin{lemma}
                                                   \label{lemma 6.4.4}
Lemma \ref{lemma 6.3.1} holds true if we replace
the restriction $\Re \lambda\geq\lambda_{0}$ in it
with $\lambda \in\Gamma$, where $\Gamma=\{\Re \lambda\geq\lambda_{0}\}
\cup\{\varepsilon_{0}|\Im \lambda|\geq-\Re \lambda+\mu_{0}\}$,
with $\varepsilon_{0}>0$ and $\mu_{0}>0$ which
depend only on $ d,\delta,p,r_{0}  $.

\end{lemma}

In the rest of the section we impose the following.
\begin{assumption}[$ p,r_{0}$] 
                                               \label{assumption 2.20.1} 
We have 
 $a^{\#}_{r_{0}}\leq\theta_{0}(d,\delta,p)$, where $\theta_{0}$
is taken in a way to accommodate  Lemmas \ref{lemma 6.3.1} and
 \ref{lemma 6.4.4}.
\end{assumption}

\begin{remark}
                                           \label{remark 6.30.1}
It is well known that if $a_{x}\in L_{d}$, then 
$a^{\#}_{r_{0}}\to0$ as $r_{0}\downarrow 0$. Therefore,
Assumption \ref{assumption 2.20.1} is weaker than
Assumption \ref{assumption 3.1.2}.
\end{remark}

On the basis of Lemma \ref{lemma 6.4.4} we can repeat
what was done in \cite{Kr_20_1} and obtain
the first part of the following result about the full operator $L$.

\begin{theorem}
                                                  \label{theorem 6.4.1}
Let $p\in(1,d)$. Then under   Assumptions  \ref{assumption 3.1.1}
and \ref{assumption 2.20.1} there   
exist
  $\lambda_{0}\geq 1, N_{0}$,
depending only on $d,\delta,p,r_{0}$, and $\nu_{b}$ (introduced below), 
  such that, for any $u\in W^{2}_{p}$
and  $\lambda \in\Gamma$,
\begin{equation}
                                 \label{6.4.6}
\sum_{r,s=1}^{d}\|D_{rs}u\|_{L_{p} } +|\lambda|\| u\|_{L_{p}}
\leq N  \|L u -\lambda u\|_{L_{p}},
\end{equation}
where $\nu_{b}$ is defined by the condition
$$
N_{1}\|bI_{|b|\geq \nu_{b}}\|_{L_{d}}\leq 1
$$
with 
a constant $N_{1}=
N_{1}(d,\delta,r_{0},p)$.
Furthermore, for any $\lambda \in\Gamma$ and $f\in L_{p}$
there is a unique  $u\in W^{2}_{p}$ such that
$\lambda u-Lu=f$.
\end{theorem}

The ``existence'' part of this theorem, as usual,
is proved by the method of continuity.

Denote by $R_{\lambda}f$ the solution $u$ from Theorem
\ref{theorem 6.4.1}. Then the fact that the norm of
$R_{\lambda}$ as an operator in $L_{p}$ decreases
as $N/|\lambda|$ for $\lambda \in\Gamma$ allows us
to use the well-known construction introduced by Hille
(\cite{Hi_48}).
We use the following facts which the reader can find, 
  for instance, in \cite{Pa_83}. For complex
$t$  in the sector $S:=\{|\Im t|<\varepsilon_{0}\Re t\}$ 
with $\varepsilon_{0}$ from Lemma \ref{lemma 6.4.4}
set
\begin{equation}
                                         \label{6.8.2}
\hat T_{t}=\frac{1}{2\pi i}\int_{\partial \Gamma}
e^{tz}R_{z}\,dz,
\end{equation}
where the integral is taken in a counter clockwise direction.
Below in this section 
$$
p\in(1,d).
$$ 

\begin{theorem}
                                                 \label{theorem 6.15.1}
(i) Formula \eqref{6.8.2} defines $\hat T_{t}$ in $S$ as
an analytic  semigroup of bounded operators
in $L_{p}$ with norms bounded 
by a constant, depending only on $\varepsilon,d,\delta,p,r_{0}$, and $\nu_{b}$, 
as long as
$t\in\{|\Im t|\leq \varepsilon \Re t, |t|\leq(\varepsilon_{0}-
\varepsilon)^{-1}\}$ for any given
$\varepsilon<\varepsilon_{0}$;

(ii) The infinitesimal generator of this semigroup is
$L$ with domain
$W^{2}_{p}$;

(iii) For $g\in W^{2}_{p}$ the
function  $\hat T_{t}g(x)$ is a unique solution of the problem
$$
\partial_{t}u(t,x)=Lu(t,x),\quad t>0,\quad \lim_{t\downarrow 0}
\|u(t,\cdot)-g\|_{L_{p}}=0
$$
 in the class of $u$ such that
$u(t,\cdot)\in W^{2}_{p}$ and (strong $L_{p}$-derivative)
$\partial_{t}u (t,\cdot) \in L_{p}$
 for each $t>0$;

(iv) For any $T\in(0,\infty)$ there is a constant $N$,
depending only on $T,d$, $\delta,p,r_{0}$, and $\nu_{b}$,  such that for each
$t\in(0,T]$ and
$f\in L_{p}$
\begin{equation}
                                         \label{6.15.1}
\|\hat T_{t}f\|_{W^{2}_{p}}\leq \frac{N}{t}
\|f\|_{L_{p}},\quad 
\|D\hat T_{t}f\|_{L_{p}}\leq \frac{N}{\sqrt{t}}
\|f\|_{L_{p}}.
\end{equation}

\end{theorem}

Actually, the second estimate in \eqref{6.15.1} is not to be found explicitly
in \cite{Pa_83} but it follows by interpolation from the first one
and the fact that $\|\hat T_{t}f\|_{L_{p}}\leq N\|f\|_{L_{p}}$.

We will also need a stability result before which we make the following.

 \begin{remark}
                                          \label{remark 6.15.1}
Let $d>p>d/2$ and $f\in W^{2}_{p}$. Then $\hat T_{t}f\in
W^{2}_{p}$ and by embedding theorems
$\hat T_{t}f$ has a modification that is  bounded and continuous in $x$, which we
still call $\hat T_{t}f$. Also $(\lambda-L)\hat T_{t}f=\hat T_{t}(\lambda-L)f\to
(\lambda-L)\hat T_{s}f$ in $L_{p} $ as $t\to s$. Hence, $\hat T_{t}f\to
\hat T_{s}f$ in $W^{2}_{p}$, which by embedding theorems
implies that $\hat T_{t}f(x)\to \hat T_{s}f(x)$ uniformly on $\bR^{d}$.
Therefore, $\hat T_{t}f(x)$ is a bounded continuous function
on $[0,T]\times\bR^{d}$ for any $T\in(0,\infty)$.

Moreover, by embedding theorems, if $u\in W^{2}_{p}$, then for any $x\in\bR^{d}$
$$
|u(x)|\leq N\big(\| u_{xx}\|_{L_{p}}+\|u\|_{L_{p}}),
$$
where $N=N(d,p)$. By substituting here $u(cx)$ in place of $u(x)$
and taking minimum of the right-hand side with respect to $c>0$
we come to the well-known estimate
$$
|u(x)|\leq N\|u_{xx}\|_{L_{p}}^{d/(2p)}\|u\|_{L_{p}}^{1-/(2p)}.
$$
Applying this and \eqref{6.15.1} yields that for $t\leq T$
and any $x\in\bR^{d}$
\begin{equation}
                                             \label{6.27.6}
|\hat T_{t}f(x)|\leq \frac{N}{t^{d/(2p)}}\|f\|_{L_{p}},
\end{equation}
where $N$ depends only on $T,d$, $\delta,p,r_{0}$, and $\nu_{b}$.
 
\end{remark}

\begin{theorem}
                                                    \label{theorem 6.23.1}

Let $d>p>d/2$ and let $a_{n},b_{n}$, $n=1,2,...$, have the same meaning as $a,b$,
respectively. Suppose that, for each $n$,  they satisfy
Assumptions \ref{assumption 3.1.1}
(with the same $\delta,\|b\|$) and \ref{assumption 2.20.1} $ (p,r_{0})$
(with the same $\theta_{0}$). Assume that $a_{n}\to a$ (a.e.)
and $b_{n}\to b$ in $L_{d}$ as $n\to\infty$. Denote by $\hat T^{n}_{t}$
the semigroups constructed on the basis of \eqref{6.8.2}
when $R_{z}$ is replaced with $R^{n}_{z}$ that is the inverse
operator to $z-L_{n}$, where $L_{n}=(1/2)a^{ij}_{n}D_{ij}+b^{i}_{n}D_{i}$.
Then for any $t>0$ and $f\in L_{p}$ we have
$\hat T^{n}_{t}f\to\hat T_{t}f$ in $W^{2}_{p}$ and, hence,
uniformly on $\bR^{d}$ as $n\to\infty$.
\end{theorem}

Proof. It suffices to prove the convergence in $W^{2}_{p}$
and formula \eqref{6.8.2} and estimate \eqref{6.4.6}
and the dominated convergence theorem show that it suffices to prove
that $\|R^{n}_{z}f-R_{z}f\|_{W^{2}_{p}}\to 0$  for
$z\in \Gamma$. In light of \eqref{6.4.6}
$$
\|R^{n}_{z}f-R_{z}f\|_{W^{2}_{p}}\leq N\|(z-L^{n})(
R^{n}_{z}f-R_{z}f)\|_{L_{p}}=N\|(L-L^{n}) 
 R_{z}f \|_{L_{p}}
$$
$$
\leq N\|\,|a^{n}-a|\,|(R_{z}f)_{xx}|\|_{L_{p}}+N
\|\,|b^{n}-b|\,|(R_{z}f)_{x }|\|_{L_{p}},
$$
where the constants $N$ are independent of $n$. In the last sum
the first term tends to zero by
the dominated convergence theorem. Concerning the second one,
observe that by the H\"older and Sobolev inequalities
$$
\|\,|b^{n}-b|\,|(R_{z}f)_{x }|\|_{L_{p}}\leq\|b^{n}-b\|_{L_{d}}
\|(R_{z}f)_{x }\|_{L_{pd/(d-p)}}
$$
$$
\leq N(d,p)
\|b^{n}-b\|_{L_{d}}
\|(R_{z}f)_{xx }\|_{L_{p }}.
$$
Therefore, it also goes to zero and the theorem is proved.

\mysection{Relation of $\hat T_{t}$ to a diffusion process}
                                                  \label{section 6.18.1}

Fix $p\in[ d_{0},d)$, where $d_{0}=d_{0}(d,\delta)\in(d/2,d)$
is taken from   \cite{Kr_19_1}, and
  suppose that Assumptions \ref{assumption 3.1.1}
and \ref{assumption 2.20.1} $(p,r_{0})$ are satisfied.

Define  $\Omega=C([0,\infty),\bR^{d})$
and for $\omega=\omega_{\cdot}\in \Omega$ define
$x_{t}(\omega)=\omega_{t}$. Also set
$\cM_{t}=\cN_{t}=\sigma(x_{s}:s\leq t)$.
Let $X=(x_{t},\infty,\cM_{t},P_{x})$ be a Markov process
corresponding to the operator $L$ constructed in \cite{Kr_19_1}
(we need only Assumptions \ref{assumption 3.1.1} for that).
We know from \cite{Kr_19_1} that, for each $x_{0}\in\bR^{d}$, 
with $P_{x_{0}}$-probability one 
\begin{equation}
                                                 \label{11.29.2}
x _{t}=x_{0}  +\int_{0}^{t}\sqrt{a (x_{s})}\,dB_{s}
+\int_{0}^{t}b (x_{s}) \,ds,
\end{equation}
where $B_{t}$ is a $d$-dimensional Wiener process on $\Omega$
relative to $\cN^{x_{0}}_{t}$ that are the completions of $\cN_{t}$
with respect to $P_{x_{0}}$.

\begin{lemma}
                                                       \label{lemma 6.30.1}
There is an extension of the  
 probability space $(\Omega, N^{x_{0}}_{\infty},P_{x_{0}})$ that carries a $d_{1}$-dimensional
Wiener process $w_{t}$   such that the above $x_{t}$
satisfies  \eqref{6.15.2}.
\end{lemma}

Proof. Enlarge the probability space $(\Omega, N^{x_{0}}_{\infty},P_{x_{0}})$
in such a way that it will carry a $d_{1}$-dimensional Wiener
process $\hat B_{t}$ the first $d$ coordinate of which coincide
with those of $B_{t}$. Then introduce $\hat \sigma^{k}=a^{-1/2}\sigma^{k}$
and observe that, for each $x$, the vectors $\xi_{i}(x)=(\hat\sigma^{i1}(x),...,
\hat\sigma^{id_{1}}(x))$, $i=1,...,d$, are orthogonal to each other
 and have unit length.
By using the Gram-Schmidt procedure it is not hard to complement them in such a way
that
$\xi_{i}(x)$, $i=1,...,d_{1}$, are orthogonal to each other,
   have unit length,
and are Borel with respect to $x$. In that case the matrix
$Q(x)$, having as rows the $\xi_{i}(x)$'s, is orthogonal and
$$
w_{t}:=\int_{0}^{t}Q^{*}(x_{s})\,d\hat B_{s}
$$
is a $d_{1}$-dimensional Wiener process. 
After that it only
remains to note that  $\hat \sigma^{k} Q^{rk}=e_{r}I_{r\leq d}$,
where $e_{r}$ is the $r$th basis vector, so that
$$
\sigma^{k}(x_{s})\,dw^{k}_{s}=a^{1/2} (x_{s})\hat
\sigma^{k}(x_{s})Q^{rk}(x_{s})\,d\hat B^{r}_{s}= a^{1/2} (x_{s})\,dB_{s}.
$$
The lemma is proved.

We know from \cite{Kr_20_1} that under Assumption \ref{assumption 2.20.1} $(p,r_{0})$  
solutions of \eqref{6.15.2} are weakly unique and therefore talking about the
properties of  solutions of \eqref{6.15.2} we may use some results
from \cite{Kr_20} about
 the process $X$.

The following result regarding $X$ is taken from
\cite{Kr_20}.
\begin{lemma}
                                                 \label{lemma 6.16.1}
Denote
$$
   T_{t}f(x)=E_{x}f(x_{t}).
$$
Then (Theorem 4.8 of \cite{Kr_20}) for any $q\geq d_{0}$ there are constants $N$ and $\mu>0$,
depending only on $d,q,\delta$, and $\|b\|$, such that
for any Borel nonnegative $f$ given 
on $\bR^{d}$ and $t>0$ we have
\begin{equation}
                                         \label{12.2.1} 
    T_{t}f(0)\leq Nt^{-d/(2q)} 
\|\Phi_{t}f\|_{L_{q}},
\end{equation}
where $\Phi_{t}(x)=\exp(-\mu |x|/\sqrt{t})$.  Furthermore
(Corollary 4.9 of \cite{Kr_20}), for   $q\geq d_{0}$ such that $q>d/2+1$
there exists a constant $N=N(q,d,\delta,\|b\|)$
such that for any  $T\in(0,\infty)$ and nonnegative Borel $f(t,x)$
given on $[0,T]\times\bR^{d} $ we have
\begin{equation}
                                                \label{12.6.1}
 E_{0}\int_{0}^{T}f(t,x_{t})\,dt\leq NT^{(q-1)/q-d/(2q)}
\| \Phi _{T}f\|_{L_{q}([0,T]\times\bR^{d})}.
\end{equation}
Finally (Lemma 6.4 of \cite{Kr_20} and \eqref{12.2.1}), $T_{t}f(x)$ is a continuous
(even locally H\"older continuous) function of $(t,x)\in(0,\infty)
\times\bR^{d}$ if $f\in L_{q}$ with $q\geq d_{0}$.

\end{lemma}
 
Note that if $u\in W^{1,2}_{q}([0,T]\times\bR^{d})$
and $q>d/2+1$, then $u$ has a modification which
is bounded and continuous on $[0,T]\times\bR^{d}$.
Therefore, talking about $u$ of class $W^{1,2}_{q}([0,T]\times\bR^{d})$
we will always mean this modification. 

\begin{theorem}[It\^o's formula]
                                             \label{theorem 6.7.1}
Let $q\geq d_{0}$ and $q>d/2+1$, and let 
$u\in W^{1,2}_{q}([0,T]\times\bR^{d})$. Then  
with  probability one for all $t\in[0,T]$ we have
\begin{equation}
                                                     \label{6.16.2}
u(t,x_{t})=u(0,x_{0})+\int_{0}^{t}(\partial_{t}+L) u(s,x_{s})\,ds
+\int_{0}^{t}\sigma^{ik}D_{i}u(s,x_{s}) \,dw^{k}_{s},
\end{equation}
where the stochastic integral is a square integrable
martingale on $[0,T]$ (and $x_{t}$ is a solution of  \eqref{6.15.2}).
\end{theorem}

This theorem is proved by using \eqref{12.6.1}
in the same way as Theorem 1.3 of \cite{Kr_19_1}
is proved on the basis of Theorem 2.6 of \cite{Kr_19_1}.

Recall that $p\in [d_{0},d)$
and $d_{0}\in(d/2,d)$, so that there are values of $p>d/2+1$
since $d\geq 3$.
 
\begin{theorem} 
                                             \label{theorem 6.16.2}
Let   $p>d/2+1$, $T\in(0,\infty)$, and 
$f\in L_{p} \cap L_{2p}$. Then 

(i) For each  $t>0$ and $x\in\bR^{d}$, we have 
$\hat T_{t}f(x)=   T_{t}f(x)$;

(ii) For each  $t>0$, for solutions of
  \eqref{6.15.2}, 
with  probability one   we
have
\begin{equation}
                                                     \label{6.16.3}
 f(x_{t})=   T_{t}f(x_{0})+\int_{0}^{t}\sigma^{ik}D_{i}
   T_{t-s}f(x_{s})\,dw^{k}_{s},
\end{equation}
where $\sigma^{ik}D_{i}
   T_{t-s}f(x )=\big(\sigma^{ik}D_{i}
   T_{t-s}f\big)(x )$ and similar notation is also used below;

(iii) For each  $t>0$ and $x\in\bR^{d}$
\begin{equation}
                                                     \label{6.17.1}
   T_{t}f^{2}(x)=(   T_{t}f(x))^{2}+
\sum_{k}\int_{0}^{t}   T_{s}\Big[\Big(\sum_{i}
\sigma^{ik}D_{i}   T_{t-s}f\Big)^{2}\Big](x)\,ds.
\end{equation}

\end{theorem}

Proof. If $f\in W^{2}_{p}$, then $u(s,x):=\hat T_{t-s}f(x)$, $s\leq t$,
satisfies the condition of Theorem \ref{theorem 6.7.1}
and we get \eqref{6.16.3} with $\hat T_{\cdot}$ in place of $   T_{\cdot}$
by that theorem.
By taking the expectations of both sides we get 
that $T_{t}f(x_{0})=  \hat  T_{t}f(x_{0})$. This holds for any $x_{0}$
and yields 
\eqref{6.16.3} as is. By taking the expectations of
the squares of both sides of \eqref{6.16.3} we obtain
\eqref{6.17.1}. Thus, all assertions of the theorem
are true if $f\in W^{2}_{p}$.

Assertion (i) holds for any  $f\in L_{p}$,
which is seen from the fact that  according to
embedding theorems and \eqref{12.2.1}
both $\hat T_{t}f(x)$ and $   T_{t}f(x)$
are bounded linear functionals on $L_{p}$
and $ W^{2}_{p}$ is dense in $L_{p}$.

Then, as $f^{n}\in  W^{2}_{p}$ tend to $f$
in $L_{p}\cap L_{2p}$, $   T_{t-s}f^{n}
\to    T_{t-s}f$ in $W^{2}_{p}$ for $s<t$ (see
\eqref{6.15.1}). By embedding theorems ($p\geq d/2$)
$D   T_{t-s}f^{n}
\to D   T_{t-s}f$ in $L_{2p}$ and in light of \eqref{12.2.1}
$$
   T_{s}\Big[\Big(\sum_{i}
\sigma^{ik}D_{i}   T_{t-s}f^{n}\Big)^{2}\Big](x)
\to    T_{s}\Big[\Big(\sum_{i}
\sigma^{ik}D_{i}   T_{t-s}f\Big)^{2}\Big](x)
$$
for any $0<s<t$ and $x\in\bR^{d}$. 
Furthermore, $(f^{n})^{2}\to f^{2}$ in $L_{p}$ and, due to \eqref{12.2.1},
$T_{t}(f^{n})^{2}(x)\to T_{t}f^{2}(x)$.
It follows by Fatou's lemma
(and \eqref{6.17.1}) that
\begin{equation}
                                                     \label{6.17.2}
   T_{t}f^{2}(x)\geq (   T_{t}f(x))^{2}+
\sum_{k}\int_{0}^{t}   T_{s}\Big[\Big(\sum_{i}
\sigma^{ik}D_{i}   T_{t-s}f\Big)^{2}\Big](x)\,ds.
\end{equation}

Hence, the right-hand side of \eqref{6.16.3}
is well defined. Furthermore,
$$
E \Big|\int_{0}^{t}\sigma^{ik}D_{i}
   T_{t-s}f(x_{s})\,dw^{k}_{s}-
\int_{0}^{t}\sigma^{ik}D_{i}
   T_{t-s}f^{n}(x_{s})\,dw^{k}_{s}\Big|^{2}
$$
$$
=\sum_{k}\int_{0}^{t}   T_{s}\Big[\Big(\sum_{i}
\sigma^{ik}D_{i}   T_{t-s}(f-f^{n})\Big)^{2}\Big](x_{0})\,ds
$$
$$
\leq    T_{t}(f-f^{n})^{2}(x_{0})- (   T_{t}(f-f^{n})(x_{0}))^{2},
$$
where the inequality is due to \eqref{6.17.2}. The last expression  
tends to zero, which allows us to get \eqref{6.16.3}
by passing to the limit in its version with $f^{n}$
in place of $f$. After that \eqref{6.17.1} follows as above.
The theorem is proved.

\begin{remark}
                                                     \label{remark 6.28.1}
In light of Theorem \ref{theorem 6.16.2} (i)
estimate \eqref{12.2.1} is weaker in what concerns the
restriction on $q$ than \eqref{6.27.6}.
However, \eqref{12.2.1} is proved for just
measurable $\sigma^{k}$. 
\end{remark}

\mysection{A criterion for strong solutions 
  of
It\^o's equations}
                                                \label{section 7.3.1}

In this section 
$$
p\in (d_{0},d),\quad p>d/2+1
$$ 
and we suppose that Assumptions \ref{assumption 3.1.1}
and \ref{assumption 2.20.1}  ($ p,r_{0}$) are satisfied.

Recall the setting from the beginning of the article.
We are given   a complete probability space $(\Omega,\cF,P)$ 
with an increasing filtration of 
$\sigma$-fields $\cF_{t}\subset \cF$, that are complete.  
We are also given  a $d_{1}$-dimensional Wiener process
$w_{t}$  relative to
$\{\cF_{t}\}$. Finally, for an $x_{0}\in \bR^{d}$
we are given a solution
  $x_{t}$
   of \eqref{6.15.2}. We know from Lemma \ref{lemma 6.30.1}
that such a situation is quite realistic and we also know that
  the solution
is weakly unique. In particular, it has the same distributions
as the process $x_{t}$ from Section \ref{section 6.18.1}
 has relative to $P_{x_{0}}$

For further discussion we need the following,
in which $\cP$ is the $\sigma$-field of predictable sets
and $\cB(0,\infty)$ is the Borel $\sigma$-field in $(0,\infty)$.

\begin{lemma}
                                              \label{lemma 6.17.1}
 Assume that
for $s,r\in (0,\infty)$, $\omega\in\Omega$ we are given
a real-valued function $g(s,r)=g(s,r,\omega)$, $s\in(0,\infty)$, $(r,\omega)
\in(0,\infty)\times\Omega$ which is measurable
with respect to  $\cB(0,\infty)
\otimes \cP$ and such that
 for each $s$
$$
E\int_{0}^{\infty}g^{2}(s,r)\,dr<\infty.
$$
Then there is a function $m_{s,t}=m(s,t,\omega)$ on $[0,\infty)\times\big(
[0,\infty) \times\Omega\big)$
measurable with respect to $\cB(0,\infty)\otimes \cP$,
continuous in $t$ for each $(s,\omega)$ and such that
  for each $s$ it is martingale starting from zero and, moreover,
for each $s$
(a.s.) for all $t\geq0$
\begin{equation}
                                                     \label{6.17.4}
m_{s,t}=\int_{0}^{t}g(s,r)\,dw_{r}.
\end{equation}
\end{lemma}

Proof. Introduce
$$
\Omega_{s}=\{\omega:\int_{0}^{\infty}g^{2}(s,r)\,dr<\infty\},\quad 
\hat g(s,r)=I_{\Omega_{s}}g(s,r),
$$
$$
  B_{t}(s)=\int_{0}^{t}\hat{g}^{2}(s,r)\,dr.
$$
Observe that
 $P(\Omega_{s})=1$ so that $\Omega_{s}\in\cF_{0}$.
Also $B_{\infty}(s)<\infty$ for any $s$ and $\omega$.

By Lemma 2.6 of \cite{Kr_11} there exists
a function $m_{s,t}$ on $[0,\infty)^{2}\times\Omega$ with the
properties described in the statement of the lemma
but satisfying \eqref{6.17.4} with $\hat g$ in place of $g$.
Since $P(\Omega_{s})=1$ the integrals of $\hat g$ and $g$
coincide with probability one and the lemma is proved.

\begin{remark}
                                             \label{remark 6.17.1}
As we have noted if $f\in L_{p}$, then,
for any $t>0$, we have  $T_{t}f\in W^{2}_{p}$, and hence ($p>d/2$),
$DT_{t}f\in L_{p}\cap L_{2p}$. Therefore, we can apply
Theorem \ref{theorem 6.16.2} and write that for any $s<t$
(a.s.)
$$
\sigma^{ik}D_{i}
   T_{t-s}f(x_{s})=T_{s}(\sigma^{ik}D_{i}
   T_{t-s}f)(x_{0})
$$
\begin{equation}
                                                     \label{6.17.20}
+\int_{0}^{s}\sigma^{jm}D_{j}T_{s-r}\big(
\sigma^{ik}D_{i}
   T_{t-s}f\big)(x_{r})\,dw^{m}_{r}.
\end{equation}
After that we want to substitute the result into \eqref{6.16.3}
to get
$$
f(x_{t})=T_{t}f(x_{0}) 
+\int_{0}^{t}T_{s}(\sigma^{ik}D_{i}
   T_{t-s}f)(x_{0})\,dw^{k}_{s}
$$
\begin{equation}
                                                     \label{6.17.3}
+\int_{0}^{t}\Big(\int_{0}^{s}\sigma^{jm}D_{j}T_{s-r}\big(
\sigma^{ik}D_{i}
   T_{t-s}f\big)(x_{r})\,dw^{m}_{r}\Big)dw^{k}_{s}.
\end{equation}
The formal objection to do that is that we should know
that the integral in \eqref{6.17.20} is, for instance,
predictable  as a function of $(\omega,s)$ and this may not
happen if we allow any version of the stochastic integral
to be taken for each $s$. However,  set $h^{k}(s,x)=
I_{s<t}\sigma^{ik}D_{i}
   T_{t-s}f(x)$ and
  consider
\begin{equation}
                                                     \label{6.17.5}
I^{k}(s,u)=\int_{0}^{u}I_{r<s }\sigma^{jm}D_{j}T_{s-r}
h^{k}(s,\cdot)(x_{r})\,dw^{m}_{r}.
\end{equation}
This is the sum over $ m$ of stochastic integrals and
$$
E \int_{0}^{\infty}I_{r<s }
\Big|\sum_{j}\sigma^{jm} D_{j}T_{s-r}h^{k}(s,\cdot)(x_{r})\Big|^{2}\,dr
$$
$$
=E \int_{0}^{s} 
\Big|\sum_{j}\sigma^{jm} D_{j}T_{s-r}h^{k}(s,\cdot)(x_{r})\Big|^{2}\,dr
\leq T_{s}\Big( \big(
h^{k}(s,\cdot)\big)^{2}\Big)(x_{0}),
$$
where the inequality is due to \eqref{6.17.1}.
It follows from Lemma \ref{lemma 6.17.1} that $I(s,u)=I(s,u,\omega)$
has a version which we denote again $I(s,u)$,
that is continuous in $u$ for each $s,\omega$ and
measurable with respect to $\cB (0,\infty) \otimes \cP $.
Then $I^{k}(s,s)$ is predictable and we take this modification
of the right-hand side of \eqref{6.17.5} in the right-hand side
of \eqref{6.17.3} thus justifying \eqref{6.17.3}.

Then we want to repeat this procedure. Introduce  
\begin{equation}
                                                         \label{6.26.5}
Q^{k}_{t}f(x)=\sigma^{ik}(x)D_{i}T_{t}f(x).
\end{equation}
In this notation \eqref{6.16.3} and \eqref{6.17.3} become, respectively,
$$
 f(x_{t})=  
T_{t}f(x_{0})+\int_{0}^{t}Q^{k_{1}}_{t-t_{1}}f(x_{t_{1}})\,dw^{k_{1}}_{t_{1}};
$$
$$
f(x_{t})=T_{t}f(x_{0}) 
+\int_{0}^{t}T_{t_{1}}Q^{k_{1}}_{t-t_{1}}(x_{0})\,dw^{k_{1}}_{t_{1}}
$$
$$
+\int_{0}^{t}\Big(\int_{0}^{t_{1}}Q^{k_{2}}_{t_{1}-t_{2}} 
Q^{k_{1}}_{t-t_{1}}f (x_{t_{2}})\,dw^{k_{2}}_{t_{2}}\Big)dw^{k_{1}}_{t_{1}}.
$$
By induction we obtain that for any $n\geq1$ (a.s.) for all $t\geq0$ ($t_{0}=t$)
$$
f(x_{t})=T_{t}f(x_{0}) 
+\sum_{m=1}^{n}\int_{t>t_{1}>...>t_{m}}T_{t_{m}}
Q^{k_{m}}_{t_{m-1}-t_{m}}\cdot...\cdot
Q^{k_{1}}_{t-t_{1}}f(x_{0})\,dw^{k_{m}}_{t_{m}}\cdot...\cdot dw^{k_{1}}_{t_{1}}
$$
\begin{equation}
                                                                \label{6.18.1}
+\int_{t>t_{1}>...>t_{n+1}}
Q^{k_{n+1}}_{t_{n}-t_{n+1}}\cdot...\cdot
Q^{k_{1}}_{t-t_{1}}f(x_{t_{n+1}})\,dw^{k_{n+1}}_{t_{n+1}}\cdot...\cdot dw^{k_{1}}_{t_{1}},
\end{equation}
where by the expressions like 
$$
\int_{t>t_{1}>...>t_{m}}:::\,dw^{k_{m}}_{t_{m}}\cdot...\cdot dw^{k_{1}}_{t_{1}}
$$
we mean
$$
\int_{ 0}^{t}\,dw^{k_{1}}_{t_{1}}\int_{  0}^{t_{1}}\,dw^{k_{2}}_{t_{2}}
...\int_{ 0}^{t_{m-1}}:::\,dw^{k_{m}}_{t_{m}}.
$$
By taking expectations of the squares of the sides in
\eqref{6.18.1} we conclude that
$$
T_{t}f^{2}(x_{0})=\big(T_{t}f(x_{0})\big)^{2}
$$
$$
+
\sum_{m=1}^{n}\int_{t>t_{1}>...>t_{m}}\big[T_{t_{m}}
Q^{k_{m}}_{t_{m-1}-t_{m}}\cdot...\cdot
Q^{k_{1}}_{t-t_{1}}f(x_{0})\big]^{2}\,d t_{m} \cdot...\cdot d t_{1} 
$$
\begin{equation}
                                                                \label{6.27.7}
+\int_{t  >t_{1}>...>t_{n+1 } }\sum_{k_{1},...,k_{n +1}}T_{t_{n+1}}
\big[
Q^{k_{n }}_{t_{n }-t_{n +1}}\cdot...\cdot
Q^{k_{1}}_{t -t_{1}}f\big]^{2}(x_{0} )\,d t_{n+1 } \cdot...\cdot d t_{1}.
\end{equation}

In particular,
 the sequence of
$$
\int_{t  >t_{1}>...>t_{n } }\sum_{k_{1},...,k_{n }}T_{t_{n}}
\big[
Q^{k_{n }}_{t_{n-1}-t_{n }}\cdot...\cdot
Q^{k_{1}}_{t -t_{1}}f\big]^{2}(x_{0})\,d t_{n } \cdot...\cdot d t_{1}
$$
is decreasing.

\end{remark}

\begin{remark}
                                                          \label{remark 6.27.1}
It turns out that proving {\em directly\/} that
each term in the right-hand side of \eqref{6.27.7}
is finite presents significant difficulties. However,
observe that, due to \eqref{6.15.1}
and \eqref{12.2.1}, for $p\in(d_{0},d)$ and $f\in L_{p} $  
we have
$$
 \big|T_{t_{m}}
Q^{k_{m}}_{t_{m-1}-t_{m}}\cdot...\cdot
Q^{k_{1}}_{t-t_{1}}f(x)\big| \leq\frac{N}{t_{m}^{d/(2p)}(t_{m-1}-t_{m})^{1/2}
\cdot...\cdot(t-t_{1})^{1/2}}\|f\|_{L_{p}}, 
$$
where $N$ depends only on $m$, $d$, $\delta$,   $\|b\|$, and $\nu_{b}$.
Furthermore,
$$
\int_{t>t_{1}>...>t_{m}}\frac{1}{t_{m}^{d/(2p)}(t_{m-1}-t_{m})^{1/2}
\cdot...\cdot(t-t_{1})^{1/2}}\,d t_{m} \cdot...\cdot d t_{1} <\infty.
$$
 
\end{remark}

Recall that
 $\cF^{w}_{t}$ is the completion of $\sigma(w_{s}:s\leq t)$. Remark \ref{remark 6.17.1}
allows us to repeat literally some arguments in \cite{VK_76}  
and leads to the following results.

\begin{theorem}
                                             \label{theorem 6.18.1}
Let $f\in L_{p}\cap L_{2p}$, $t>0$. Then
$$
E \big(f(x_{t})\mid\cF^{w}_{t}\big)=T_{t}f(x_{0})
$$
$$
+\sum_{m=1}^{\infty}\int_{t>t_{1}>...>t_{m}}T_{t_{m}}
Q^{k_{m}}_{t_{m-1}-t_{m}}\cdot...\cdot
Q^{k_{1}}_{t-t_{1}}f(x_{0})\,dw^{k_{m}}_{t_{m}}\cdot...\cdot dw^{k_{1}}_{t_{1}},
$$
where the series converges in the mean square sense.
\end{theorem}

For $n\geq 1$, $t>0$, and $s_{1},...,s_{n}>0$  
define
\begin{equation}
                                                                \label{6.26.1}
Q_{s_{n},...,s_{1}}f(x)=\sum_{k_{1},...,k_{n }}
\big[
Q^{k_{n }}_{s_{n }}\cdot...\cdot
Q^{k_{1}}_{s_{1}}f\big]^{2}(x ).
\end{equation}

\begin{theorem}
                                             \label{theorem 6.18.2}
Let $f\in L_{p}\cap L_{2p}$, $t_{0}>0$. Then
$f(x_{t_{0}})$ is $\cF^{w}_{t_{0}}$-measurable iff
\begin{equation}
                                                                \label{6.18.10}
\lim_{n\to\infty}  \int_{t_{0}>t_{1}>...>t_{n } }T_{t_{n }}
Q_{t_{n-1}-t_{n},...,t_{0}-t_{1}}f(x_{0})
\,d t_{n } \cdot...\cdot d t_{1}=0.
\end{equation}
Furthermore, under either of the above equivalent conditions
$$
 f(x_{t}) =T_{t}f(x_{0})
$$
\begin{equation}
                                                                \label{7.9.1}
+\sum_{m=1}^{\infty}\int_{t>t_{1}>...>t_{m}}T_{t_{m}}
Q^{k_{m}}_{t_{m-1}-t_{m}}\cdot...\cdot
Q^{k_{1}}_{t-t_{1}}f(x_{0})\,dw^{k_{m}}_{t_{m}}\cdot...\cdot dw^{k_{1}}_{t_{1}}.
\end{equation}

\end{theorem}

\begin{theorem}
                                             \label{theorem 6.18.3}
If equation \eqref{6.15.2} has two solutions which
are not indistinguishable, then it does not have any
strong solution. In particular, if \eqref{6.15.2} has at least
one strong solution, then the solution is unique.

\end{theorem}

\begin{theorem}
                                             \label{theorem 7.1.1}
If equation \eqref{6.15.2} has a strong
solution on one probability space then it has a strong solution
on any other probability space carrying a $d_{1}$-dimensional
Wiener process.

\end{theorem}

\begin{remark}
                                                        \label{remark 6.26.1}
By making the change of variables $t_{k}=s_{k}+...+s_{n}$,
$k=1,...,n$ we rewrite \eqref{6.18.10} as
\begin{equation}
                                                                \label{6.26.3}
\lim_{n\to\infty}  \int_{S_{n}(t_{0}) }T_{s_{n }}
Q_{s_{n-1} ,...,s_{1},t_{0}-(s_{1}+...+s_{n})}f(x_{0})
\,d s_{n } \cdot...\cdot d s_{1}=0,
\end{equation}
where $S_{n}(t_{0})=\{(s_{1},...,s_{n}):s_{k}\geq0,s_{1}+...+s_{n}< t_{0}\}$

\end{remark}

The sequence under the limit sign in \eqref{6.26.3},
call it $u_{n}(t_{0})$, is decreasing
for any $t_{0}$ (and $x$). Therefore, its limit will be zero
for almost any $t$ if  
$$
\lim_{n\to\infty} \int_{0}^{\infty}u_{n}(t_{0})e^{-\nu t_{0}}\,dt_{0}=0,
$$
where $\nu>0$ is any number.  In that case, actually, the limit of $u_{n}(t_{0})$
is zero for all $t_{0}$, since, in light of Theorem \ref{theorem 6.18.2},
$f(x_{t_{0}})$ is $\cF^{w}_{t_{0}}$-measurable for almost
all $t_{0}$, and by continuity, for all $t_{0}$.
In this way after simple change of variables we come to the following.

\begin{theorem}
                                           \label{theorem 6.26.1}
Let $f\in L_{p}\cap L_{2p}$. Then
$f(x_{t })$ is $\cF^{w}_{t }$-measurable for any $t>0$
if there exists a $\nu>0$ such that
\begin{equation}
                                                                \label{6.26.4}
\int_{0}^{\infty}e^{-\nu s_{n}}T_{s_{n }} 
\Big( \int_{R^{n}_{+} }e^{-\nu (s_{n-1}+...+s_{0})}
Q_{s_{n-1} ,...,s_{0} }f
\,d s_{n-1 } \cdot...\cdot d s_{0}\Big)(x_{0})\,ds_{n}\to 0
\end{equation}
as $n\to \infty$, where $R^{n}_{+}=(0,\infty)^{n }$, which in light of \eqref{12.2.1} holds for any $x_{0}$ if
\begin{equation}
                                                                \label{6.26.50}
  \Big\| \int_{R^{n }_{+} }e^{-\nu (s_{n-1}+...+s_{0})}
Q_{s_{n-1} ,...,s_{0} }f(x)
\,d s_{n-1 } \cdot...\cdot d s_{0}\Big\|^{p}_{L_{p}} \to 0.
\end{equation}
\end{theorem}

We are going to prove that \eqref{6.26.50} holds under
Assumptions \ref{assumption 3.1.1} and \ref{assumption 3.1.2}
by showing that the series composed of the left-hand sides
of \eqref{6.26.50} converges.

\begin{remark}
                                               \label{remark 7.10.1}
The criterion \eqref{6.18.10} is 
proved under Assumptions \ref{assumption 3.1.1}
and \ref{assumption 2.20.1} $(p,r_{0})$, assumptions,
 which involve the $\sigma^{k}$'s
only implicitly and it turns out that for some choice
of the $\sigma^{k}$'s \eqref{6.18.10} may hold and for another fail
to hold. To illustrate this we take $b\equiv 0$. In that 
case the restriction $p<d$ disappears along with $d\geq3$
(which is a consequence of $p<d$ and $p>d/2+1$). Then we take $d_{1}=d= 2$
and following \cite{KZ_81} set $\sigma^{1}(x)=x/|x|$, $\sigma^{2}(x)=x^{*}/|x|$,
where $x^{*}=(-x^{2},x^{1} )$
for $x\ne 0$, $\sigma^{ik}(0)=\delta^{ik}$. Then $a^{ij}(x)=\delta^{ij}$,
equation \eqref{6.15.2} has a solution for any $x_{0}$
(see, for instance, Lemma \ref{lemma 6.30.1}), and each solution
is a Wiener process starting from $x_{0}$. For $x_{0}\ne0$
the solutions are strong and, hence, \eqref{6.18.10}
holds, because the solution never reaches the origin, the point where
$\sigma^{k}$ are not smooth. However, for $x_{0}=0$
there are no strong solutions, because, as is easy to see,
rotation in $x^{1}x^{2}$ coordinates  by any angle brings 
any solution
it to another solution
of the same equation. Therefore, for $x_{0}=0$ equation 
\eqref{6.18.10} does not hold. 

Also observe that in this example $\sigma^{k}\in W^{1}_{d-\varepsilon,\loc}$
for any $\varepsilon\in(0,1)$ but not for $\varepsilon=0$.
One can construct similar examples for $d\geq3$ starting
from the following   with $d=3$, $d_{1}=9$, and $\sigma^{k}$'s
that are the $k$th columns of the matrix
$$
\frac{1}{|x|}
\begin{pmatrix}
x^{1} &  x^{2} & x^{3} & 0 &  0   & 0   & 0  & 0  & 0 \\
0 & 0 & 0 & x^{1}  & x^{2}   & x^{3}   & 0  & 0  & 0 \\
0 & 0 & 0 & 0  & 0   & 0  & x^{1}& x^{2} & x^{3} 
\end{pmatrix}.
$$
Again $a^{ij}=\delta^{ij}$, $\sigma^{k}\in W^{1}_{d-\varepsilon,\loc}$
for any $\varepsilon\in(0,1)$ but not for $\varepsilon=0$,
and, if $x_{t}$ is a solution of \eqref{6.15.2} with $x_{0}=0$,
then $-x_{t}$ is also a solution of \eqref{6.15.2} with $x_{0}=0$.
\end{remark}

\mysection{Some estimates in the case of $C^{\infty}$ coefficients}
                                                \label{section 7.3.2}

We suppose that $\sigma^{k},b$ satisfy Assumption
\ref{assumption 3.1.1} and are infinitely differentiable
with each derivative bounded.

Let $(\Omega,\cF,P)$ be a complete probability space,
let $\{\cF_{t}\}$ be an increasing filtration of 
$\sigma$-fields $\cF_{t}\subset \cF$, that are complete.
Let 
$w_{t}$ be a $d_{1}$-dimensional Wiener process relative to
$\{\cF_{t}\}$. We also assume that there is a $(d+1)$ independent $d$-dimensional
Wiener, relative to  $\{\cF_{t}\}$, process $ B^{(0)}_{t},...,B^{(d)}_{t}$ 
independent
of $w_{t}$. Take $x,\eta\in \bR^{d}$, a nonnegative bounded
infinitely differentiable  $K_{0}$ with each derivative bounded
given on $\bR^{d}$,
and consider the following system   
\begin{equation}
                                                        \label{6.20.3}
x_{t}=x+\int_{0}^{t}\sigma^{k}(x_{s})\,dw^{k}_{s}+
\int_{0}^{t}b(x_{s})\,ds,
\end{equation}
$$
\eta_{t}=\eta+\int_{0}^{t}\sigma^{k}_{(\eta_{s})}(x_{s})\,dw^{k}_{s}
+\int_{0}^{t}b_{(\eta_{s})}(x_{s})\,ds
$$
\begin{equation}
                                                        \label{6.20.4}
+\int_{0}^{t}K_{0}(x_{s})\,dB^{(0)}_{s}+\int_{0}^{t}K_{0}(x_{s}) \eta^{k}_{s} 
\,dB^{(k)}_{s} .
\end{equation}
As is well known, \eqref{6.20.3} 
 has a unique solution which we denote by $x_{t}(x)$.
By substituting it into \eqref{6.20.4} we see that the coefficients
of \eqref{6.20.4} grow linearly in $\eta$ and hence 
\eqref{6.20.4} also has a unique solution which we denote by
$\eta_{t}(x,\eta)$. By the way observe that equation \eqref{6.20.4}
is linear with respect to $\eta_{t}$. Therefore
$\eta_{t}(x,\eta)$ is an affine function of $\eta$.
For the uniformity of notation we set $x_{t}(x,\eta)=x_{t}(x)$.
It is also well known 
  (see, for instance, Sections 2.7 and 2.8 of
\cite{Kr_77}) that, as a function of $x$ and $(x,\eta)$, the processes
$x_{t}(x)$ and $\eta_{t}(x,\eta)$
 are infinitely differentiable in an appropriate sense
(specified below),
their derivatives satisfy the equations which are obtained by formal
differentiation of \eqref{6.20.3} and \eqref{6.20.4},
respectively, and, for any $n\geq0,T\in(0,\infty)$,  $l_{k},\xi_{k}\in\bR^{d}$,
$k=1,...,n$ (if $n\geq1$),
$x,\eta\in \bR^{d}$, and $q\geq 1$,
\begin{equation}
                                                        \label{6.21.1}
 E\sup_{t\leq T}\Big|\Big(\prod_{k=1}^{n}(lb)D _{(l_{k},\xi_{k})}
\Big)(x_{t},\eta_{t})(x,\eta)
\Big|^{q}\leq N(1+|\eta|^{m}),
\end{equation}
where $N $ is a certain constant  independent of $(x,\eta)$,
$m=m(n,q)$, and,
for instance,
by $(lb)D _{(l,\xi)} \eta_{t} (x,\eta)$ we mean a process $\zeta_{t}$
such that, for any $q\geq1$ and $S\in(0,\infty)$
$$
\lim_{\varepsilon\downarrow0}
E\sup_{t\leq S}\big|\zeta_{t}-\varepsilon^{-1}
\big(\eta_{t} (x+\varepsilon l,\eta+\varepsilon\xi)
-\eta_{t} (x,\eta)\big)\big|^{q}=0.
$$
\begin{lemma}
                                                     \label{lemma 6.21.10}
Let $f(x,\eta)$ be infinitely differentiable and such that
each of its derivatives grows  as $|x|+|\eta|\to\infty$
not faster than polynomially. Then the function
$u(t,x,\eta):=Ef\big((x_{t},\eta_{t})(x,\eta)\big)$
is infinitely differentiable in $(x,\eta)$ and 
each of its derivatives is continuous in $t$ and is by absolute
value bounded on each finite time interval 
by a constant times $(1+|x|+|\eta|)^{m}$
for some $m$. Furthermore, $u(t,x,\eta)$ is continuously
differentiable in $t$ and for $t\geq0$ and $(x,\eta)\in\bR^{2d}$
$$
 \partial_{t}u(t,x,\eta)= (1/2)\sigma^{ik}\sigma^{jk}(x)u_{x^{i}x^{j}} (t,x,\eta)
+\sigma^{ik}\sigma_{(\eta)}^{jk}(x)u_{x^{i}\eta^{j}} (t,x,\eta)
$$
$$
+(1/2)\sigma_{(\eta)}^{ik} \sigma_{(\eta)}^{jk}(x)u_{\eta^{i}\eta^{j}}(t,x,\eta)
+(1/2)K^{2}_{0}(x)(1+|\eta|^{2})\delta^{ij}u_{\eta^{i}\eta^{j}}(t,x,\eta)
$$
\begin{equation}
                                                        \label{6.21.3}
+b^{i}(x)u_{x^{i}} (t,x,\eta)+b^{i}_{(\eta)}(x)u_{\eta^{i}} (t,x,\eta)
=:\check L(x,\eta)u(t,x,\eta).
\end{equation}
\end{lemma}

The first assertion of this lemma follows easily from what is said before it.
Then the fact that \eqref{6.21.3} holds follows from the Markov
property of $(x_{t},\eta_{t})$ and from the first assertion.
The claimed property of $\partial_{t}u$ follows from \eqref{6.21.3}.

\begin{lemma}
                                                     \label{lemma 6.21.1}
 Let $\eta\in\bR^{d}$   and
 $\xi_{t}(x,\eta)=(lb)D_{\eta}x_{t}(x)$. Then

(i) $\xi_{t}(x,\eta)$ satisfies \eqref{6.20.4}
with $K_{0}\equiv0$.

(ii) For any $t$
\begin{equation}
                                                        \label{6.21.4}
 \xi_{t}(x,\eta)=E\big(\eta_{t}(x,\eta)\mid\cF^{w}_{t}\big).
\end{equation}

(iii) If $f(x)$ is infinitely differentiable with bounded derivatives,
then
\begin{equation}
                                                        \label{6.21.5}
Ef_{(\eta_{t}(x,\eta))}(x_{t}(x))\Big(=
E\big(f_{(\eta_{t}(x,\eta))}\big)(x_{t}(x))\Big)
=\big(Ef(x_{t}(x))\big)_{(\eta)}.
\end{equation}
\end{lemma}

Proof. Assertion (i) is well known (see, for instance, \cite{Kr_77}).
The right-hand side of \eqref{6.21.4} satisfies \eqref{6.20.4}
with $K_{0}\equiv0$   owing to the linearity of $g_{(\eta)}$
in $\eta$ and independence of $B_{\cdot}$ and $w_{\cdot}$. 
Therefore, due to uniqueness,
assertion (ii) follows from (i). Assertion (iii) follows from (ii)
and the fact that (see, for instance, \cite{Kr_77})
$$
\big(Ef(x_{t}(x))\big)_{(\eta)}=Ef_{(\xi_{t}(x,\eta))}(x_{t}(x)).
$$
The lemma is proved.

Now follows one of the most important computations. The idea behind it
is the following. If we formally differentiate both parts of \eqref{7.9.1}
in the direction $\eta$ and then take the expectations of the squares
of both sides, then we obtain an equality in \eqref{6.21.6} below
if we also replace on the left $\eta_{t}(x,\eta)$ by $\xi_{t}(x,\eta)$.
Then the inequality follows from Lemma \ref{lemma 6.21.1}.

\begin{lemma}
                                                     \label{lemma 6.21.01}
 Let $x,\eta\in\bR^{d}$ and let $f (x)$
be infinitely differentiable with bounded derivatives.
 Then for any $t\in(0,\infty)$
$(t_{0}=t$)
$$
E\big[f_{(\eta_{t}(x,\eta))}(x_{t}(x))\big]^{2}
\geq\Big[(T_{t}f(x))_{(\eta)}\Big]^{2}
$$
\begin{equation}
                                                        \label{6.21.6}
+\sum_{n=1}^{\infty}\sum_{k_{1},...,k_{n}}
\int_{t>t_{1}>...>t_{n}}\Big[\big(T_{t_{n}}Q^{k_{n}}_{t_{n-1}-t_{n}}
\cdot...\cdot  Q^{k_{1}}_{t-t_{1}}f(x)\big)_{(\eta)}\Big]^{2}\,dt_{n}
\cdot...\cdot dt_{1}.
\end{equation}
\end{lemma}

Proof. Introduce the notation $\check T_{t}u(x,\eta)=Eu\big((x_{t},\eta_{t}) (x,\eta)\big)$.
Then, similarly to \eqref{6.16.3}, for smooth bounded
$u(x,\eta)$  by dropping for simplicity the arguments $x$ and $\eta$  
in $x_{s}(x)$ and $\eta_{s}(x,\eta)$, we get
$$
u(x_{t},\eta_{t})=\check T_{t}u(x,\eta)
+\int_{0}^{t}K_{0}D_{\eta^{i}}\check T_{t-t_{1}}u(x_{t_{1}},\eta_{t_{1}})\,
 \big( dB^{i(0)}_{t_{1}}+\eta^{k}_{t_{1}}dB^{i(k)}_{t_{1}}\big)
$$
$$
+\int_{0}^{t}\Big[\sigma^{ik}(x_{t_{1}})D_{x^{i}}
\check T_{t-t_{1}}u(x_{t_{1}},\eta_{t_{1}})+
\sigma^{ik}_{(\eta_{t_{1}})}(x_{t_{1}})D_{\eta^{i}}
\check T_{t-t_{1}}u(x _{t_{1} },\eta_{t_{1}} )\Big]\,dw^{k}_{t_{1}}.
$$
It follows that
$$
Eu^{2}(x_{t},\eta_{t})\geq \big(\check T_{t}u(x,\eta)\big)^{2}
$$
\begin{equation}
                                                        \label{6.21.7}
+\sum_{k}\int_{0}^{t}E\Big[\sigma^{ik}(x_{t_{1}})D_{x^{i}}
\check T_{t-t_{1}}u(x_{t_{1}},\eta_{t_{1}})+
\sigma^{ik}_{(\eta_{t_{1}})}(x_{t_{1}})D_{\eta^{i}}
\check T_{t-t_{1}}u(x _{t_{1}},\eta_{t_{1}})\Big]^{2}\,dt_{1}.
\end{equation}
By using Fatou's lemma and  estimates like \eqref{6.21.1}
one easily carries \eqref{6.21.7} over to smooth $u(x,\eta)$
whose derivatives have no more than polynomial growth
as $|x|+|\eta|\to\infty$. In particular, one can apply
\eqref{6.21.7} to $u(x,\eta)=f_{(\eta)}(x)$. Then,
after noting that in light of \eqref{6.21.5} in that case
$$
\sigma^{ik}(x )D_{x^{i}}
\check T_{t-t_{1}}u(x ,\eta )+
\sigma^{ik}_{(\eta)}(x )D_{\eta^{i}}
\check T_{t-t_{1}}u(x   ,\eta  )
$$
$$
=\sigma^{ik}(x )D_{x^{i}}(T_{t-t_{1}}f(x))_{(\eta)}
+\sigma^{ik}_{(\eta )}(x )D_{\eta^{i}}(T_{t-t_{1}}f(x))_{(\eta)}
$$
$$
=\big(\sigma^{ik}(x )D_{x^{i}} T_{t-t_{1}}f(x)\big)_{(\eta)}
=\big(Q^{k}_{t-t_{1}}f(x)\big)_{(\eta)},
$$
we obtain
$$
E\big[f_{(\eta_{t} )}(x_{t} )\big]^{2}
\geq\Big[(T_{t}f(x))_{(\eta)}\Big]^{2}
+\sum_{k_{1}}
\int_{0}^{t}E\big[(Q^{k_{1}}_{t-t_{1}}f)_{(\eta_{t} )}
(x_{t_{1}} )\big]^{2}\,dt_{1}.
$$
By applying this formula to $Q^{k_{1}}_{t-t_{1}}f$ in place of $f$
we get
$$
E\big[f_{(\eta_{t})}(x_{t})\big]^{2}
\geq\Big[(T_{t}f(x))_{(\eta)}\Big]^{2}+\sum_{k_{1}}
\int_{0}^{t}E\big[(T_{t_{1}}Q ^{k_{1}}_{t-t_{1}}f (x  ))_{(\eta)}\big]^{2}\,dt_{1} 
$$
$$
+\sum_{k_{1},k_{2}}\int_{0}^{t}dt_{1}
\int_{0}^{t_{1}}E
\big[(Q^{k_{2}}_{t_{1}-t_{2}}Q^{k_{1}}_{t-t_{1}}f)_{(\eta_{t_{2})}}(x_{t_{2}})
\big]^{2}\,dt_{2}.
$$

 Using induction yields that for any $n\geq1$ $(t_{0}=t$)
$$
E\big[f_{(\eta_{t})}(x_{t})\big]^{2}
\geq\Big[(T_{t}u(x))_{(\eta)}\Big]^{2}
$$
$$
+\sum_{m=1}^{n}\sum_{k_{1},...,k_{m}}
\int_{t>t_{1}>...>t_{m}}\Big[\big(T_{t_{m}}Q^{k_{m}}_{t_{m-1}-t_{m}}
\cdot...\cdot  Q^{k_{1}}_{t-t_{1}}f(x)\big)_{(\eta)}\Big]^{2}\,dt_{m}
\cdot...\cdot dt_{1}
$$
$$
+\sum_{k_{1},...,k_{n+1}}
\int_{t>t_{1}>...>t_{n+1}}E\Big[\big( Q^{k_{n+1}}_{t_{n }-t_{n+1}}
\cdot...\cdot  Q^{k_{1}}_{t-t_{1}}f\big)_{(\eta_{t_{n+1}})}
(x_{t_{n+1}})\Big]^{2}\,dt_{n+1}
\cdot...\cdot dt_{1}.
$$
This yields \eqref{6.21.6} and proves the lemma.

Next, we want to estimate the left-hand side of \eqref{6.21.6}
which according to Lemma \ref{lemma 6.21.10} satisfies
\eqref{6.21.3}. 

In the future we will be interested in estimating not only
 the left-hand side of \eqref{6.21.6} but a slightly more
general quantity. Therefore, we take an infinitely
differentiable $f(x,\eta)\geq0$ such that
for an $m>0$ and a constant $N$
$$
\big(|f|+|f_{x}|+|f_{\eta}|+|f_{xx}|+|f_{x\eta}|+
|f_{\eta\eta}|\big)(x,\eta)\leq N(1+|\eta|)^{m}
$$
for all $x,\eta$  and such that $f(x,\eta)=0$
for all $\eta$ if $|x|\geq R$ for some $R>0$. Then denote 
$u(t,x,\eta)=\check T_{t}f(x,\eta)$. According to what was said before
Lemma \ref{lemma 6.21.1} and in that lemma, if we denote
$(l_{t},\xi_{t})=(lb)D _{(l,\xi)}(x_{t}, \eta_{t}) (x,\eta)$,
then 
$$
u_{(l,\xi)}(t,x,\eta)=Ef_{(l_{t},\xi_{t})}(x_{t},\eta_{t})(x,\eta).
$$
In particular, it follows that
$$
|u_{(l,\xi)}(t,x,\eta)|\leq P^{1/2}(|x_{t}(x)|\leq R) 
\Big(E |f_{(l_{t},\xi_{t})}(x_{t},\eta_{t})(x,\eta)|^{2}\Big)^{1/2}.
$$
Here the second factor on the right is estimated by using \eqref{6.21.1}.
The first factor is less than $2\exp(-\mu\text{dist}^{2}(x,B_{R})/t)$ by Theorem 2.10
of \cite{Kr_19}, where $\mu>0$ depends only on $d,\delta$, and $\|b\|$.
Similar estimates are available for the second-order derivatives of
$u (t,x,\eta)$. 
More precisely,
observe that  
there exist  constants  $ \mu >0$, $\kappa=\kappa(m)$, and a function
$M(t)$ bounded on each time interval $[0,T]$ such that for all $t $,
  $x,\eta$  we have
$$
 |u(t,x,\eta)|+| u_{x}(t,x,\eta)|+
| u_{\eta}(t,x,\eta)|
$$
\begin{equation}
                                                        \label{6.21.8}
+| u_{xx}(t,x,\eta)|
+| u_{x\eta}(t,x,\eta)| +| u_{\eta\eta}(t,x,\eta)| 
\leq M(t) e^{-\mu |x| } (1+|\eta|^{2})^{\kappa}.
\end{equation}
This justifies the integration we perform below.

Introduce
$$
h  = (1+| \eta|^{2})^{-\kappa-d} 
$$ 
and observe that for a constant $N=N(d,\kappa)$ we have
$$
|\eta|\,|h _{\eta}|\leq N \kappa h ,\quad 
(1+| \eta|^{2})|
h^{\kappa}_{\eta \eta }|\leq  h .
$$

\begin{theorem}
                                                   \label{theorem 6.21.1}
Let   $q\geq2$ and suppose that the above $u\geq0$.
 Then  there is a constant $N_{0}
\geq1$,
depending only on    $d$, $d_{1}$, $\delta$,  $m$, and $q$,
such that for any $\lambda\geq 1$ satisfying
\begin{equation}
                                                        \label{7.1.1}
N_{0}\Big(\sum_{k}\| \sigma^{k}_{x}I_{| \sigma^{k}_{x}|>\lambda}\|_{L_{d}}
+\|bI_{|b|>\lambda}\|_{L_{d}}\Big)\leq 1
\end{equation}
there exists  a constant  $N$,
depending only on   $\lambda $, $\| \sigma^{k}_{x}\|$, $d$, $\delta$,  $m$,  and $q$,
  and there is a function $K_{0}$ 
such that for $t\geq0$
\begin{equation}
                                                        \label{6.21.90}
\int_{\bR^{2d}}h ( \eta) u^{q}(t,x,\eta) \,dxd\eta\leq 
e^{Nt}\int_{\bR^{2d}}h ( \eta)f^{q}( x,\eta)\,dxd\eta.
\end{equation}

\end{theorem}

The proof of this theorem proceeds as usual
 by multiplying \eqref{6.21.3} by $h (\eta)u^{q-1}(t,x,\eta)$ 
and integrating by parts
over $[0,t]\times\bR^{2d}$. The integral of the left-hand side
is
$$
q^{-1}\int_{\bR^{2d}}h ( \eta)u^{q}(t,x,\eta)\,dxd\eta-q^{-1}
\int_{\bR^{2d}}h ( \eta)f^{q}( x,\eta)\,dxd\eta.
$$
Therefore, in light of Gronwall's inequality, to prove the theorem
it suffices to prove the following estimate.  

\begin{lemma}
                                                        \label{lemma 6.21.5}
Let  $q\geq2$ and $\kappa\geq0$. 
Then  there is a constant $N_{0}\geq 1$,
depending only on    $d$, $d_{1}$, $\delta$,   $\kappa$, and $q$,
such that for any $\lambda\geq 1$ satisfying
\begin{equation}
                                                        \label{7.1.10}
N_{0}\Big(\sum_{k}\| \sigma^{k}_{x}I_{| \sigma^{k}_{x}|>\lambda}\|_{L_{d}}
+\|bI_{|b|>\lambda}\|_{L_{d}}\Big)\leq 1
\end{equation}
there exists a  constant  $N$,
depending only on   $\lambda $, $\|\sigma^{k}_{x}\|$, 
$d$, $d_{1}$, $\delta$,  $\kappa$, and $q$,
  and there is a function $K_{0}$  
  such that  for any a smooth 
function $v(x,\eta)\geq0$ (independent of $t$), for which condition
 \eqref{6.21.8} is satisfied with $v$ in place of $u$ and some $M$,
we have
\begin{equation}
                                                        \label{6.21.9}
\int_{\bR^{2d}}h ( \eta) v^{q-1} ( x,\eta)\check Lv(x,\eta) \,dxd\eta
\leq N\int_{\bR^{2d}}h ( \eta) v^{q}( x,\eta) \,dxd\eta .
\end{equation}

\end{lemma}

Proof. For simplicity of notation we drop the arguments $x,\eta$. 
We also write $U\sim V$ if
their integrals over $\bR^{2d}$ coincide, and $U\prec V$ if the integral of
$U$ is less than or equal to that of $V$.  Below the constants called $N$,
sometimes with indices, depend
 only on $d,d_{1}$, $\delta$,  $\kappa$, and~$q$ unless specifically noted otherwise.

Set   $w=v^{q/2}$ and note simple formulas:
$$
v^{q-1}v_{x}=(2/q)ww_{x},\quad v^{q-2}v_{x^{i}}v_{x^{j}}
=(4/q^{2})w_{x^{i}}w_{x^{j}}.
$$
Then
denote by $\check L_{1}$ the sum of the first-order terms in $\check L$
 and observe that integrating by parts shows that
$$
h v^{q-1}  b^{i}_{(\eta)} v_{\eta^{i}}  \sim -(1/q)
h _{  \eta^{i}}  b^{i}_{(\eta)}v^{q}
-(1/q)h b^{i}_{x^{i}} v^{q}
$$
$$
\sim   (2/q)\eta^{k}h _{\eta^{i} } b^{i}w  w_{x^{k}} 
+(2/q)h b^{i}ww_{x^{i}} .
$$
Hence,
$$
hv^{q-1}\check L_{1}v \sim  
 (2/q)\eta^{k}h_{\eta^{i} } b^{i}w  w_{x^{k}} 
+(4/q)hb^{i}ww_{x^{i}} .
$$
We take a number $\lambda\geq1 $ and write $b=\hat b+\check b$,
where $\hat b=bI_{|b |>\lambda}$. Following \cite{St_65}
we observe that by the H\"older and Sobolev inequalities ($d\geq3$)
$$
\int_{\bR^{d}}|\hat b^{i}ww_{x^{k}}|\,dx \leq
\Big(\int_{\bR^{d}} |w_{x} |^{2}\,dx\Big)^{1/2}
\Big(\int_{\bR^{d}} |\hat b|^{2}|w|^{2}\,dx\Big)^{1/2}
$$
\begin{equation}
                                                        \label{6.22.2}
\leq  
\Big(\int_{\bR^{d}} |w_{x} |^{2}\,dx\Big)^{1/2}
 \|\hat b\|_{L_{d}} \|w\|_{L_{2d/(d-2)}} \leq
N\|\hat b\|_{L_{d}}\int_{\bR^{d}} |w_{x} |^{2}\,dx ,
\end{equation}
where  $N$ depends only on $d$. Since $|\eta|\,|h_{\eta}|\leq
N(\kappa,d)h$,
it follows that 
$$
\eta^{k}h_{\eta^{i} }\hat b^{i}ww_{x^{k}}
\prec N\|\hat b\|_{L_{d}}h|w_{x} |^{2}.
$$ 
Similarly,
$(4/q)h\hat b^{i}ww_{x^{i}}\prec N\|\hat b\|_{L_{d}} h|w_{x} |^{2}$.
We estimate the remaining terms in $h v^{q-1}\check L_{1}v$ roughly  like
$$
|\check  b^{i}ww_{x^{k}}|\leq \lambda|w|\,|w_{x} |
\leq \varepsilon|w_{x} |^{2}+\varepsilon ^{-1}\lambda^{2}|w|^{2}
$$
and conclude that, for any $\varepsilon>0$,
\begin{equation}
                                                        \label{6.22.1}
 h^{q}v^{q-1}\check L_{1}v\prec (N\|\hat b\|_{L_{d}}+
\varepsilon)h |w_{x} |^{2}+N\varepsilon ^{-1}\lambda^{2}h|w|^{2}.
\end{equation}

Starting to deal with the second order derivatives note that
$$
h v^{q-1}(1/2)\sigma^{ik}\sigma^{jk} v_{x^{i}x^{j}} \sim -
((q-1)/2)v^{q-2}h\sigma^{ik}v_{x^{i}} \sigma^{jk} v_{x^{j} } 
$$
$$
-(1/2)h\big[\sigma^{ik}_{x^{i}}\sigma^{jk}+
\sigma^{ik}\sigma^{jk}_{x^{i}}\big]v^{q-1}v_{x^{j}} =
-((2q-2)/q^{2}) h\sigma^{ik}w_{x^{i}} \sigma^{jk} w_{x^{j} } 
$$
$$
-(1/q)h\big[\sigma^{ik}_{x^{i}}\sigma^{jk}+
\sigma^{ik}\sigma^{jk}_{x^{i}}\big]ww_{x^{j}}  
\leq -(1/q) h\sigma^{ik}w_{x^{i}} \sigma^{jk} w_{x^{j} }
$$
$$
+h\Big|\big[\sigma^{ik}_{x^{i}}\sigma^{jk}+
\sigma^{ik}\sigma^{jk}_{x^{i}}\big]ww_{x^{j}}\Big|,
$$
where the inequality (to simplify the writing) is due to the fact that $q\geq2$.
In this inequality the first term on the right is dominated in the sense
of $\prec$ by 
$$
-(1/q)\delta h|w_{x}|^{2}
$$
 (see \eqref{7.18.1}).
The remaining term contains $ww_{x^{i}}$ and we treat it like above writing
$\sigma^{k}_{x}=\hat \sigma^{k}+\check \sigma^{k}$,
where $\hat \sigma^{k}=\sigma^{k}_{x}I_{| \sigma^{k}_{x} |> \lambda}$.
Then we get
$$
h v^{q-1}(1/2)\sigma^{ik}\sigma^{jk} v_{x^{i}x^{j}}\prec
N\lambda^{2}\varepsilon^{-1}h|w|^{2}
$$
\begin{equation}
                                                        \label{6.22.3}
-\Big[(1/q)\delta-N\sum_{k}\|\hat \sigma^{k}\|_{L_{d}}-\varepsilon\Big] h|w_{x}|^{2}.
\end{equation}

Next,
$$
h v^{q-1}\sigma^{ik}\sigma_{(\eta)}^{jk} v_{x^{i}\eta^{j}}  \sim
-(q-1)h\sigma^{ik}v^{q-2} v_{\eta^{j}} \sigma_{(\eta)}^{jk} v_{x^{i}}
$$
$$
-v^{q-1} v_{x^{i}} \big[h_{\eta^{j}}  \sigma^{ik}\sigma_{(\eta)}^{jk}
+h \sigma^{ik}\sigma_{x^{j}}^{jk}]=-((4q-4)/q^{2})
h\sigma^{ik}  w_{\eta^{j}} \sigma_{(\eta)}^{jk} w_{x^{i}}
$$
$$
-(2/q)ww_{x^{i}} \big[h_{\eta^{j}}  \sigma^{ik}\sigma_{(\eta)}^{jk}
+h\sigma^{ik}\sigma_{x^{j}}^{jk}].
$$
We estimate the first term on the right roughly using
$$
|\sigma^{ik} w_{\eta^{j}} \sigma_{(\eta)}^{jk} w_{x^{i}}|\leq
\varepsilon|w_{x}|^{2}+N\varepsilon^{-1}|\eta|\sum_{k}|\sigma^{k}_{x}|^{2}
|w_{\eta}|^{2}.
$$
The second term  contains $ww_{x^{i}}$ and allows the same handling as before.
 Therefore
$$
h v^{q-1}\sigma^{ik}\sigma_{(\eta)}^{jk} v_{x^{i}\eta^{j}} \prec
 \Big(\varepsilon+N\sum_{k}\|\hat \sigma^{k}\|_{L_{d}}\Big) h|w_{x}|^{2}
$$
\begin{equation}
                                                        \label{6.22.4}
+N\varepsilon^{-1}h|\eta|\sum_{k}|\sigma^{k}_{x}|^{2}|w_{\eta}|^{2}
+N\lambda^{2}\varepsilon^{-1}h|w|^{2}.
\end{equation}
 
The last term in $h v^{q-1}\check  Lv$ containing $\sigma$ is
$$
h v^{q-1}(1/2)\sigma_{(\eta)}^{ik} \sigma_{(\eta)}^{jk} v_{\eta^{i}\eta^{j}} \sim
-((q-1)/2)h\sigma_{(\eta)}^{ik} v^{q-2} v_{ \eta^{j}}  \sigma_{(\eta)}^{jk}v_{\eta^{i}}
$$
$$
-(1/2)v^{q-1}\sigma_{(\eta)}^{ik} v_{\eta^{i} }\big[h_{\eta^{j}}  \sigma_{(\eta)}^{jk}
+h \sigma_{x^{j}}^{jk}\big]-(1/(2q))h (v^{q})_{\eta^{i}}
 \sigma_{x^{j}}^{ik} \sigma_{(\eta)}^{jk}
$$
$$
\prec Nh(|\eta|^{2}|w_{\eta}|^{2}+w^{2})\sum_{k}|\sigma^{k}_{x}|^{2} +I,
$$
where
$$
I=-(1/(2q))h(w^{2})_{\eta^{i}}
 \sigma_{x^{j}}^{ik} \sigma_{(\eta)}^{jk}
$$
$$
\sim
(1/(2q))w\sigma_{x^{j}}^{ik}\big[h_{\eta^{i}} \sigma_{(\eta)}^{jk}
+h\sigma_{x^{i}}^{jk}\big]\prec Nh\sum_{k}|\sigma^{k}_{x}|^{2}w^{2}.
$$
To estimate the last term we basically use the derivation
of \eqref{6.22.2}. We have
\begin{equation}
                                                          \label{7.5.2}
\int_{\bR^{d}}|\hat\sigma^{k}|^{2}w^{2}\,dx\leq
\|\hat\sigma^{k}\|^{2}_{L_{d}}\|w\|_{L_{2d/(d-2)}}^{2}\leq N
\|\hat\sigma^{k}\|^{2}_{L_{d}}
\int_{\bR^{d}}|w_{x}|^{2}\,dx.
\end{equation}    
Below we show how to choose the constant $N_{0}$ in the condition
\eqref{7.1.10} under which our assertion is true. But observe that with
any such choice $\|\hat\sigma^{k}\|^{2}_{L_{d}}\leq 1$ and therefore,
(just to keep some uniformity in our estimates) \eqref{7.5.2}
implies that
\begin{equation}
                                                          \label{7.5.5}
\int_{\bR^{d}}|\hat\sigma^{k}|^{2}w^{2}\,dx \leq N
\|\hat\sigma^{k}\| _{L_{d}}
\int_{\bR^{d}}|w_{x}|^{2}\,dx.
\end{equation} 
Hence,
$$
I\prec Nh\lambda^{2}w^{2}+Nh\sum_{k}\|\hat\sigma^{k}\| _{L_{d}}|w_{x}|^{2}
$$
and
$$
h v^{q-1}(1/2)\sigma_{(\eta)}^{ik} \sigma_{(\eta)}^{jk} v_{\eta^{i}\eta^{j}}
\prec Nh |\eta|^{2}|w_{\eta}|^{2} \sum_{k}|\sigma^{k}_{x}|^{2}
$$
\begin{equation}
                                                                 \label{7.9.3}
+Nh\lambda^{2}w^{2}+Nh\sum_{k}\|\hat\sigma^{k}\| _{L_{d}}|w_{x}|^{2}
\end{equation}

Finally,
$$
h v^{q-1}(1/2)K^{2}_{0}(1+|\eta|^{2})\delta^{ij}v_{\eta^{i}\eta^{j}} \sim
-((2q-2)/q^{2})hK^{2}_{0}(1+|\eta|^{2})|w_{\eta }|^{2}
$$
$$
-(2/q)K^{2}_{0}\big(h(1+|\eta|^{2})\big)
_{\eta^{i}}ww_{\eta^{i}}
$$
$$
\sim -((2q-2)/q^{2})hK^{2}_{0}(1+|\eta|^{2})|w_{\eta }|^{2}+(1/q) w^{2}
K_{0}^{2}\delta^{ij}\big(h(1+|\eta|^{2})\big)
_{\eta^{i}\eta^{j}}
$$
\begin{equation}
                                                                 \label{7.9.4}
\leq -(1/q)hK^{2}_{0}(1+|\eta|^{2})|w_{\eta }|^{2}+N w^{2}
K_{0}^{2}h.
\end{equation}

By combining \eqref{6.22.1}, \eqref{6.22.3},
 \eqref{6.22.4}, \eqref{7.9.3}, and \eqref{7.9.4}, 
 we see that for any $\varepsilon\in(0,1]$
$$
h v^{q-1}\check Lv\prec  N\varepsilon ^{-1}\lambda^{2}h|w|^{2}
$$
$$
-\Big[(1/q)\delta-N_{1}\Big(\sum_{k}\|\hat \sigma^{k}\|_{L_{d}}
+\|\hat b\|_{L_{d}}\Big)-N_{2}\varepsilon\Big] h|w_{x}|^{2}
$$
\begin{equation}
                                                              \label{7.5.3}
+|w_{\eta}|^{2}h(1+|\eta|^{2})\Big[N_{3}\varepsilon^{-1}\sum_{k}|\sigma^{k}_{x}|^{2}
-(1/q)K_{0}^{2}\Big]
+N_{4} w^{2}
K_{0}^{2}h.
\end{equation}
Take and fix $\varepsilon$ so that $N_{2}\varepsilon\leq \delta/(2q)$.
After that set
$$
K_{0}^{2}=1+N_{3}q\varepsilon^{-1}\sum_{k}|\sigma^{k}_{x}|^{2}
$$
(1 is added to guarantee the smoothness of $K_{0}$)
and observe that similarly to \eqref{7.5.2} and \eqref{7.5.5}
$$
N_{4} w^{2}
K_{0}^{2}h=N_{4} w^{2}
  h+Nh w^{2}\sum_{k}|\sigma^{k}_{x}|^{2}\prec N_{5}h
\sum_{k}\|\hat \sigma^{k}\| _{L_{d}}|w_{x}|^{2}+Nh\lambda^{2}w^{2}.
$$
Then \eqref{7.5.3} becomes
$$
h v^{q-1}\check Lv\prec  N \lambda^{2}h|w|^{2}
-\Big[(1/(2q))\delta-(N_{1}+N_{5})\Big(\sum_{k}\|\hat \sigma^{k}\|_{L_{d}}
+\|\hat b\|_{L_{d}}\Big) \Big] h|w_{x}|^{2}.
$$
We can certainly believe that $N_{1}\geq1$, 
take $N_{0}$ in \eqref{7.1.10} to be equal to $(2q/\delta)(N_{1}+N_{5})$
($\geq1$),
and conclude
$$
h v^{q-1}\check Lv\prec  N \lambda^{2}h|w|^{2}.
$$
The lemma is proved.

We finish the section with an approximation result.

\begin{lemma}
                                                        \label{lemma 7.2.1}
Let  Assumptions \ref{assumption 3.1.1} and
\ref{assumption 3.1.2} be satisfied. Then
there are sequences $\sigma^{k}(n)$, $b(n)$, $n=1,2,...$, $k=1,...,d_{1}$,
of infinitely differentiable functions with each derivative bounded
having the same meanings as $\sigma^{k}$, $b$ in the beginning
of the article, satisfying Assumptions \ref{assumption 3.1.1} and
\ref{assumption 3.1.2} with $\delta/2$ 
in place of $\delta $ and the same $ \|b\|$  and $ \| \sigma^{k}_{x}\|$
for sufficiently large $n$,
and such that $\sigma^{k}(n)\to \sigma^{k}$ as $n\to\infty$ (a.e.)
and $ \sigma^{k}_{x}(n), b(n)\to  \sigma^{k}_{x}, b$ in $L_{d}$ as $n\to\infty$.
\end{lemma}

Proof. Take a nonnegative $\zeta\in C^{\infty}_{0}$
with unit integral and support in $B_{1}  $ and set
$\zeta_{n}(x)=n^{d}\zeta(nx)$, $u(n,x)=u(x)*\zeta(nx)$.
Then the well-known properties of convolutions
imply all stated properties apart from what concerns
\eqref{7.18.1}. 

Denote by $\sigma$ the $d\times d_{1}$-matrix
whose columns are the $\sigma^{k}$'s and note that
$$
|\sigma^{*}(n,x)\lambda|\leq \zeta(nx)|\sigma^{*}(x)\lambda|\leq\delta^{-1/2}|\lambda|.
$$
Therefore we need only prove that for sufficiently large $n$
\begin{equation}
                                                    \label{7.2.1}
|\sigma^{*}(n,x)\lambda|\geq |\lambda|\delta^{-1/2}/\sqrt2.
\end{equation}
For any $y$ we have
$$
|\sigma^{*}(n,x)\lambda|\geq |\sigma^{*}(y)\lambda|-
 |(\sigma^{*}(n,x)-\sigma^{*}(y))\lambda|\geq |\lambda|\delta^{1/2}
-|(\sigma^{*}(n,x)-\sigma^{*}(y))\lambda| 
$$
$$
\geq |\lambda|\big(\delta^{1/2}
-|(\sigma^{*}(n,x)-\sigma^{*}(y))|\big)
$$
Furthermore,
$$
\int_{\bR^{d}}|(\sigma^{*}(n,x)-\sigma^{*}(x-y)) |\zeta_{n}(y)\,dy
$$
$$
\leq \int_{\bR^{d}}\int_{\bR^{d}}
|(\sigma^{*}( x-z)-\sigma^{*}(x+y)) |\zeta_{n}(y)\zeta_{n}(z)\,dydz
$$
$$
\leq N(d)\max\zeta^{2}\osc(\sigma,B_{1/n}(x)).
$$
The latter tends to zero uniformly with respect to $x$
since $ \sigma_{x}\in L_{d}$ (cf.~Remark \ref{remark 6.30.1}).
 This certainly proves the lemma.

\mysection{Proof of Theorem \protect\ref{theorem 6.29.1}}

                                                \label{section 7.3.4}

According to Theorems \ref{theorem 6.18.3} and \ref{theorem 7.1.1}
 it suffices to prove that
at least one of the solutions of \eqref{6.15.2} is strong. We will be dealing with the
solution from Lemma \ref{lemma 6.30.1}.

Let $f\in C^{\infty}_{0}$.  First we deal with smooth
coefficients and develop necessary estimates. By Lemma \ref{lemma 6.21.01}  
and  Theorem  \ref{theorem 6.21.1} for $t\geq0$ and $q\geq2$
we have
\begin{equation}
                                                         \label{6.23.2}
\int_{\bR^{2d}}h ( \eta) u^{q}(t,x,\eta) \,dxd\eta\leq Ne^{Nt},
\end{equation}
where (and below) 
 $N$ depends only on $f$,   $d$, $d_{1}$, $\delta$, $m=m(f)$, $q$, 
and $\lambda$ defined by \eqref{7.1.1}
and
$$
u (t,x,\eta)
$$
$$
:=\sum_{n=1}^{\infty}\sum_{k_{1},...,k_{n}}
\int_{t>t_{1}>...>t_{n}}\Big[\big(T_{t_{n}}Q^{k_{n}}_{t_{n-1}-t_{n}}
\cdot...\cdot  Q^{k_{1}}_{t-t_{1}}f(x)\big)_{(\eta)}\Big]^{2}\,dt_{n}
\cdot...\cdot dt_{1}.
$$
Obviously, $u (t,x,\eta)$ is a quadratic function
of $\eta$. Hence, \eqref{6.23.2} implies that, for any $R\in(0,\infty)$

\begin{equation}
                                                         \label{6.23.3}
\int_{\bR^{d}} \sup_{|\eta|\leq R} u^{q}(t,x,\eta) \,dx \leq Ne^{Nt}R^{2q}.
\end{equation}
Observe that in notation \eqref{6.26.5} and \eqref{6.26.1}
$$
\sum_{k }u(t,x,\sigma^{k})
=\sum_{n=1}^{\infty} 
\int_{t>t_{1}>...>t_{n}}  Q_{t_{n},t_{n-1}-t_{n},
 ...,  t-t_{1}}f(x)   \,dt_{n}
\cdot...\cdot dt_{1}
$$
$$
=\sum_{n=1}^{\infty} 
\int_{S_{n}(t)}  Q_{s_{n},s_{n-1} ,  
 ...,s_{1},  t-(s_{1}+...+s_{n}}f(x)   \,ds_{n}
\cdot...\cdot ds_{1}=:\sum_{n=1}^{\infty}I_{n}(t,x)
$$
($S_{n}(t)$ is introduced in Remark \ref{remark 6.26.1}).
Next, for $\nu>0$ by Jensen's inequality
$$
\sum_{n=1}^{\infty}\int_{\bR^{d}}\Big(\int_{0}^{\infty}
e^{-\nu t}I_{n}(t,x)\,dt\Big)^{q}\,dx
\leq \nu^{1-q}\int_{0}^{\infty}e^{-\nu t}\Big(\sum_{n=1}^{\infty}
\int_{\bR^{d}}I^{q}_{n}(t,x)\,dx\Big)dt
$$
$$
\leq \nu^{1-q}\int_{0}^{\infty}e^{-\nu t}\int_{\bR^{d}}
\Big(\sum_{k }u(t,x,\sigma^{k})\Big)^{q}\,dxdt ,
$$
which thanks to \eqref{6.23.3} implies that  for appropriate $\nu$,
depending only   on $f$,   $d$, $d_{1}$, $\delta$,   $q$, 
and $\lambda$,
\begin{equation}
                                                              \label{6.27.1}
\sum_{n=1}^{\infty}\int_{\bR^{d}}\Big(\int_{0}^{\infty}
e^{-\nu t}I_{n}(t,x)\,dt\Big)^{q}\,dx
 \leq N,
\end{equation}
where   $N$ depends  only only on $f$,   $d$, $d_{1}$, $\delta$,   $q$, 
and $\lambda$.

Estimate \eqref{6.27.1} has  been derived only for
infinitely differentiable $\sigma$ and $b$. However, using
smooth approximations (Lemma \ref{lemma 7.2.1}), Theorem \ref{theorem 6.23.1}, and Fatou's
lemma prove \eqref{6.27.1} also in our general case. Indeed, although
the constant $N$ in \eqref{6.27.1} for each approximation
depends on $\lambda$, satisfying \eqref{7.1.1} for the approximating
$\sigma_{x}^{k}$ and $b$, it can be taken   the same as long as
the approximations are sufficiently close in $L_{d}$ to the original 
$\sigma_{x}^{k}$ and $b$.

Finally, by observing that
$$
\int_{0}^{\infty}
e^{-\nu t}I_{n}(t,x)\,dt=\int_{R^{n+1}_{+}}e^{-\nu(s_{0}+...+s_{n})}
Q_{s_{n},...,s_{0}}f(x)\,ds_{n}\cdot...\cdot ds_{0},
$$
 referring to Theorem \ref{theorem 6.26.1}, and taking $q=p$, we 
conclude that $f(x_{t})$ is $\cF^{w}_{t}$-measurable.
The arbitrariness of $f$ and $t$
finishes the proof.

\mysection{Proof of Theorem \protect\ref{theorem 6.29.3}}
                                                \label{section 7.3.5}

Take a bounded smooth function $f$ with compact support.
By Theorems \ref{theorem 6.29.1} and \ref{theorem 6.18.2} for any $t$
$$
 f(x_{t}(n,x(n))) =T_{t}(n)f(x(n))
$$
\begin{equation}
                                                                \label{7.2.2}
+\sum_{m=1}^{\infty}\int_{t>t_{1}>...>t_{m}}T_{t_{m}}(n)
Q^{k_{m}}_{t_{m-1}-t_{m}}(n)\cdot...\cdot
Q^{k_{1}}_{t-t_{1}}(n)f(x(n))\,dw^{k_{m}}_{t_{m}}\cdot...\cdot dw^{k_{1}}_{t_{1}},
\end{equation}
where $T_{t}(n)$ and $Q^{k}_{t}(n)$ are the operators corresponding
to $\sigma^{k}(n)$, $b(n)$. First we prove that $E|f(x_{t}(n,x(n)))
-f(x_{t} )|^{2}\to0$ as $n\to\infty$. Since $Ef^{2}(x_{t}(n,x(n)))
\to Ef^{2}(x_{t} )$ (see Theorem \ref{theorem 6.23.1}), it suffices to prove  
that $f(x_{t}(n,x(n)))\to f(x_{t} )$ weakly in $L_{2}(\Omega,\cF^{w}_{t},P)$.
Furthermore, according to \cite{It_51} the linear combinations
of the multiple It\^o integrals of the type
$$
\int_{t>t_{1}>...>t_{m}}\phi(t_{1},...,t_{m})\,dw_{t_{m}}\cdot...\cdot dw_{t_{1}},
$$
where $m$ is arbitrary and $\phi$ is an arbitrary bounded (nonrandom) Borel function,
are dense in $L_{2}(\Omega,\cF^{w}_{t},P)$.
Therefore, it suffices to prove that for all such $m$ and $\phi$
 $$
Ef\big(x_{t}(n,x(n))\big)\int_{t>t_{1}>...>t_{m}}
\phi(t_{1},...,t_{m})\,dw_{t_{m}}\cdot...\cdot dw_{t_{1}}
$$ 
 $$
\to Ef(x_{t} )\int_{t>t_{1}>...>t_{m}}
\phi(t_{1},...,t_{m})\,dw_{t_{m}}\cdot...\cdot dw_{t_{1}}.
$$
In light of \eqref{7.2.2} this is equivalent to proving that
 $$
\int_{t>t_{1}>...>t_{m}}\phi(t_{1},...,t_{m})T_{t_{m}}(n)
Q^{k_{m}}_{t_{m-1}-t_{m}}(n)\cdot...\cdot
Q^{k_{1}}_{t-t_{1}}(n)f(x(n))\,dt_{m}\cdot...\cdot dt_{1}
$$ 
 $$
\to \int_{t>t_{1}>...>t_{m}}\phi(t_{1},...,t_{m})T_{t_{m}} 
Q^{k_{m}}_{t_{m-1}-t_{m}} \cdot...\cdot
Q^{k_{1}}_{t-t_{1}} f(x_{0})\,dt_{m}\cdot...\cdot dt_{1}.
$$
This relation is indeed true, which follows by the dominated
convergence theorem from Theorem \ref{theorem 6.23.1} 
and Remark \ref{remark 6.27.1}.

Next, observe that for any $T\in(0,\infty)$ and bounded smooth
$\bR^{d}$-valued $\tilde b$ with compact support
$$
I:=\nlimsup_{n\to\infty}E\sup_{t\leq T}\Big|\int_{0}^{t}b(n,x_{s}(n,x(n)))\,ds-
\int_{0}^{t}b( x_{s} )\,ds\Big|
$$
$$
\leq\nlimsup_{n\to\infty} E\int_{0}^{T}|b(n,x_{s}(n,x(n)))-  b(x_{s} )|\,ds
$$
$$
\leq\nlimsup_{n\to\infty} E\int_{0}^{T}|b(n,x_{s}(n,x(n)))
-\tilde b(x_{s}(n,x(n)) )|\,ds
$$
$$
+\nlimsup_{n\to\infty}  \int_{0}^{T}E|\tilde b(x_{s}(n,x(n)) )-\tilde b(x_{s} )|\,ds
$$
$$
+\nlimsup_{n\to\infty} E\int_{0}^{T}|\tilde b(x_{s} ) -b( x_{s} )|\,ds.
$$
Here the middle term vanishes by the first part of the proof.
Owing to Lemma \ref{lemma 6.16.1}
the two remaining terms are majorated by
$$
N(\lim_{n\to\infty}\|b_{n}-\tilde b\|_{L_{d}}+ \|b -\tilde b\|_{L_{d}})
=N\|b -\tilde b\|_{L_{d}},
$$
that can be made arbitrarily small by an appropriate choice of $\tilde b$.
Hence, $I=0$.

Similarly,
$$
J:=\nlimsup_{n\to\infty}
\Big(E\sup_{t\leq T}\Big|\int_{0}^{t}\sigma^{k}(n,x_{s}(n,x(n)))\,dw^{k}_{s}-
\int_{0}^{t}\sigma^{k}( x_{s} )\,dw^{k}_{s}\Big|\Big)^{2}
$$
$$
\leq N\nlimsup_{n\to\infty}
  \sum_{k}E\int_{0}^{T}|\sigma^{k}(n,x_{s}(n,x(n)))-\sigma^{k}( x_{s} )|^{2}\,ds
$$
$$
\leq N\nlimsup_{n\to\infty}
 \sum_{k}E\int_{0}^{T}|\sigma^{k}(n,x_{s}(n,x(n)))-\sigma^{k}( x_{s} )| \,ds
$$
$$
\leq N\sum_{k}\|\Phi _{T}(\sigma^{k} -\hat\sigma^{k}) \|_{L_{d}},
 $$
where $\hat\sigma^{k}$ are smooth functions with compact support.
It follows that $J=0$ and together with $I=0$ this implies that
\begin{equation}
                                                            \label{7.2.3}
\lim_{n\to\infty}E\sup_{t\leq T}| x_{t}(n,x(n)) - x_{t} |=0
\end{equation}
By Corollary 1.2 of \cite{Kr_19} for any $m\geq0$
$$
E\sup_{t\leq T}|x_{t}(n,x(n))-x(n)|^{2m}+E\sup_{t\leq T}|x_{t}  -x_{0}|^{2m}\leq  
N(m,d,\delta,\|b\|)T^{m}
$$
and this along with \eqref{7.2.3} yields the result.
The theorem is proved.  

\mysection{Proof of Theorem \protect\ref{theorem 6.29.2}}
                                                \label{section 7.3.9}

 First we assume that $\sigma^{k}$ and $b$ are infinitely differentiable
with each derivative bounded.
In that case, as  it is known since \cite{BF_61}
(see also \cite{Ku_90}) that one can define $x_{t}(x)$ in such a way
that it becomes differentiable in $x$ for all $(\omega,t)$
and the derivative of $x_{t}$ in the direction of $\eta$
satisfies the same equation as $\xi_{t}(x,\eta)$ from Lemma \ref{lemma 6.21.1},
for which \eqref{6.21.4} holds. In particular, for any even $\kappa \geq2$,
and $f$ with compact support ($(\cdot,\cdot)$ is the scalar product in $\bR^{d}$)
$$
 E\big((Df)(x_{t}(x)),\eta_{t}(x,\eta)\big)^{\kappa }\geq 
E\big((Df)(x_{t}(x)),\xi_{t}(x,\eta)\big)^{\kappa }=:v(t,x,\eta).
$$
By Theorem \ref{theorem 6.21.1}, with $q=2$ there, there is a constant $m=m(\kappa )$
such that for any $\lambda>0$ satisfying \eqref{7.1.1}
there exists a constant $N$,
depending only on   $\lambda $, $d$, $\delta$, $m $,    such that for $t\geq0$
$$
\int_{\bR^{2d}}h ( \eta) v^{2}(t,x,\eta) \,dxd\eta
$$
\begin{equation}
                                                        \label{7.2.4}
\leq 
e^{Nt}\int_{\bR^{2d}}h ( \eta)|f_{(\eta)}(x)| ^{2\kappa }\,dxd\eta
=N(d,\kappa)e^{Nt}\int_{\bR^{ d}}|f_{x}(x)|^{2\kappa }\,dx=:M_{t}.
\end{equation}
Next, for any $R\in(0,\infty)$
$$
 E\int_{B_{R}}|D\big(f(x_{t}(x))\big)|^{\kappa }\,dx 
=NE\int_{B_{R}}\int_{\bR^{ d}}\big(\big(f(x_{t}(x))_{(\eta)}\big) ^{\kappa }h 
(\eta)\,d\eta
dx
$$
$$
=N \int_{B_{R}}\int_{\bR^{ d}}v(t,x,\eta)h (\eta)\,d\eta.
$$
By using \eqref{7.2.4} and H\"older's inequality we obtain that
$$
 E\int_{B_{R}}|D\big(f(x_{t}(x))\big)|^{\kappa }\,dx \leq
N(d,\kappa)R^{d/2}M^{1/2}_{t}.
$$
By Morrey's theorem (see, for instance, Theorem 
10.2.1 of \cite{Kr_08}) this implies that for 
any  $\kappa >d$
\begin{equation}
                                                        \label{7.3.1}
E \sup_{x,y\in B_{R}}\frac{|f(x_{t}(x) )-f(x_{t}(y) )|^{\kappa }
}{|x-y|^{\kappa -d}} \leq
N(d, \kappa )R^{d/2}M^{1/2}_{t}.
\end{equation} 

Note that \eqref{7.3.1} is certainly applicable to vector-valued $f$.
Fix $\rho\geq 2R$ and take a smooth $f$ with  support
in $B_{4\rho}$  
such that $f(x)=x$ for $|x|\leq 2\rho$ and $|f_{x}|\leq 2$.
Then $M_{t}\leq N(T,d,\delta,\lambda,\kappa)\rho^{d}$ and for $x,y\in B_{R}$ and $t\leq T$
$$
E|x_{t}(x)-x_{t}(y)|^{\kappa}\leq  
N\rho^{d/2}|x-y|^{\kappa-d}  
$$
$$
+N(\kappa)E\big(|x_{t}(x) |^{\kappa}+|x_{t}(x) |^{\kappa}\big)
I_{|x_{t}(x)|+|x_{t}(y)|\geq 2\rho},
$$
where $N$ depends only on $R,T,\lambda,d,\delta$, and $\kappa$. 
We estimate the first term on the right by using H\"older's inequality and
  Theorem 2.10 of \cite{Kr_19} and find that it is dominated
by   
$$
N(\kappa)P^{1/2}(|x_{t}(x)|+|x_{t}(y)|\geq 2\rho)
\Big(\big(E |x_{t}(x) |^{2\kappa}\big)^{1/2}
+\big(E |x_{t}(x) |^{2\kappa}\big)^{1/2}\Big)
$$
$$
\leq Ne^{-\mu \rho^{2}},
$$
where $N$ depends only on $R,T,\lambda,d,\delta$, and $\kappa$
and $\mu>0$ depends only on $T, d,\delta $, and $\|b\|$. Thus, for
$\rho\geq 2R$
$$
E|x_{t}(x)-x_{t}(y)|^{\kappa}\leq Ne^{-\mu \rho^{2}}+
N\rho^{d/2}|x-y|^{\kappa-d}.
$$
By taking here $\kappa>2d$ and $\mu\rho^{2}=-\ln |x-y|^{\kappa-2d}$ we find that
\begin{equation}
                                                      \label{7.11.2}
E|x_{t}(x)-x_{t}(y)|^{\kappa}\leq  
N |x-y|^{\kappa-2d},
\end{equation}
where $N$ depends only on $R,T,\lambda,d,\delta,\|b\|$, and $\kappa$,
provided that $-\ln |x-y|^{\kappa-2d}\geq 2R/\mu$. However, if
 $-\ln |x-y|^{\kappa-2d}\leq 2R/\mu$, \eqref{7.11.2} is obvious.

Estimate \eqref{7.11.2} so far is proved only for infinitely differentiable
coefficients, but usual approximations, Theorem \ref{theorem 6.29.3},
and Fatou's lemma allow us to obtain \eqref{7.11.2} in our general case
where, naturally, by
$x_{t}(x)$ we mean the strong solution of \eqref{6.15.2} with $x_{0}=x$.

In the general case we also have by Corollary 1.2 of \cite{Kr_19}
that for any $q \geq1$
\begin{equation}
                                                        \label{7.3.2}
E| x_{t}(x) - x_{s}(x)  |^{q}  \leq N|t-s|^{q/2},
\end{equation} 
where $N=N(d,\delta,\|b\|,q )$.

Now the arbitrariness of $\kappa$ and $q$ leads to the claimed result
by a version of Kolmogorov's theorem which can be found, for instance, in
\cite{Po_09} or simply derived by using
$|x_{t}(x) - x_{s}(y)|\leq |x_{t}(x) - x_{s}(x)|+|x_{s}(x) - x_{s}(y)|$.
  The theorem is proved.

\end{document}